\documentclass[]{interact}
\usepackage[percent]{overpic}
\usepackage[numbers, sort&compress]{natbib}
\bibpunct[, ]{[}{]}{,}{n}{,}{,}% Citation support using natbib.sty
\usepackage{amsmath}
\usepackage{amscd}
\usepackage{amssymb}
\usepackage{graphicx}
\usepackage{url}
\usepackage{color}
\usepackage{listings}
\usepackage{booktabs}
\usepackage{multirow}
\usepackage{textcomp}
\usepackage{mathtools}
\usepackage{subcaption}
\usepackage{hyperref}
\usepackage[utf8]{inputenc}
\usepackage{tikz}
\usepackage{pgfplots}
\usepackage{environ}
\usepackage{algorithm}
\usepackage{enumerate}
\usepackage[noend]{algpseudocode}
\usepackage{color}
\usepackage{filecontents}

\immediate\write18{bibtex \jobname}

\definecolor{blue}{rgb}{0.2980392156862745, 0.4470588235294118, 0.6901960784313725}
\definecolor{green}{rgb}{0.3333333333333333, 0.6588235294117647, 0.40784313725490196}
\definecolor{red}{rgb}{0.7686274509803922, 0.3058823529411765, 0.3215686274509804}

\theoremstyle{plain}% Theorem-like structures provided by amsthm.sty
\newtheorem{theorem}{Theorem}[section]
\newtheorem{lemma}[theorem]{Lemma}

\theoremstyle{definition}
\newtheorem{definition}[theorem]{Definition}

\theoremstyle{remark}
\newtheorem{remark}{Remark}

\definecolor{darkblue}{rgb}{0.00,0.00,0.55}
\definecolor{black}{rgb}{0.00,0.00,0.00}
\hypersetup{
    linkcolor = darkblue,
    anchorcolor = darkblue,
    citecolor = darkblue,
    filecolor = darkblue,
    urlcolor = darkblue,
    colorlinks  = true,
}

\title{Deflation for semismooth equations\thanks{This research is supported by DFG grant no.~SU 963/1-1
		``Generalized Nash Equilibrium Problems with Partial
		Differential Operators: Theory, Algorithms, and Risk
		Aversion'', by EPSRC grants EP/K030930/1 and EP/M011151/1,
		and by the EPSRC Centre For Doctoral Training in
		Industrially Focused Mathematical Modelling
		(EP/L015803/1) in collaboration with Simula Research
		Laboratory. The authors would like to
		acknowledge useful discussions with M.~C.~Ferris, D.~Klatte and C.~Kanzow, and thank M.~C.~Ferris
		for supplying a GAMS formulation of the example considered
		in section \ref{sec:gerard}.}}

\author{
	\name{A.~N. Author\textsuperscript{a}\thanks{CONTACT A.~N. Author. Email: latex.helpdesk@tandf.co.uk} and John Smith\textsuperscript{b}}
	\affil{\textsuperscript{a}Taylor \& Francis, 4 Park Square, Milton Park, Abingdon, UK; \textsuperscript{b}Institut f\"{u}r Informatik, Albert-Ludwigs-Universit\"{a}t, Freiburg, Germany}
}

\author{
  \name{Patrick E.~Farrell\textsuperscript{a}\thanks{Email: \texttt{patrick.farrell@maths.ox.ac.uk}. Corresponding author.},
  	    Matteo Croci\textsuperscript{ab}\thanks{Email: \texttt{matteo.croci@maths.ox.ac.uk}.}
  	    and
  	    Thomas M.~Surowiec\textsuperscript{c}\thanks{Email: \texttt{surowiec@mathematik.uni-marburg.de}.}
  }
  \affil{\textsuperscript{a}Mathematical Institute, University of Oxford, Oxford, UK;
  		 \textsuperscript{b}Simula Research Laboratory, Oslo, Norway;
  		 \textsuperscript{c}Fachbereich Mathematik und Informatik, Philipps-Universität Marburg, Marburg, Germany
  }
}

\begin{document}
\maketitle

\begin{abstract}
Variational inequalities can in general support distinct solutions. In this 
paper we study an algorithm for computing distinct solutions of a variational
inequality, without varying the initial guess supplied to the solver.  The
central idea is the combination of a semismooth Newton method with a deflation
operator that eliminates known solutions from consideration. Given one root
of a semismooth residual, deflation constructs a new problem for which a semismooth Newton
method will not converge to the known root, even from the same initial
guess. This enables the discovery of other roots. We prove the
effectiveness of the deflation technique under the same assumptions that
guarantee locally superlinear convergence of a semismooth Newton method. We demonstrate its utility 
on various finite- and infinite-dimensional examples drawn from constrained optimization, game theory, economics and solid mechanics.
\end{abstract}


\begin{amscode}
    65K15 Numerical methods for variational inequalities and related problems,
    65P30 Bifurcation problems,
    65H10 Systems of equations,
    35M86 Nonlinear unilateral problems and nonlinear variational inequalities of mixed type,
    90C33 Complementarity and equilibrium problems and variational inequalities (finite dimensions).
\end{amscode}

\section{Introduction}

Variational inequalities are a fundamental class of problem that arise in many
branches of applied mathematics. The problem to be
solved is: given a real reflexive Banach space $U$, a closed convex subset $K \subset U$, and an operator $Q:
K \to U^*$ mapping to the dual space $U^*$ of $U$, find $u \in K$ such that
\begin{equation} \label{eqn:vi}
\langle Q(u), v - u \rangle \ge 0 \quad \text{for all } v \in K.
\end{equation}
This is denoted by VI$(Q, K)$.  As an elementary example, consider the problem
of minimizing a differentiable function $f: \mathbb{R} \to \mathbb{R}$ over the
closed interval $I \subset \mathbb{R}$. The necessary condition for $z \in
\mathbb{R}$ to be a (local) minimum is that $z$ satisfies VI($f', I$). More
generally, the minimizers of a general smooth nonlinear program satisfy
a variational inequality, which is related to the familiar Karush--Kuhn--Tucker conditions
under a suitable constraint qualification, see e.g.~\cite{nocedal2006}.  Variational inequalities also arise naturally in
problems of solid mechanics involving contact \cite{fichera1973,kikuchi1988}, in game theory
for the calculation of Nash equilibria \cite{harker1984}, in phase separation
with nonsmooth free energy \cite{blowey1991}, and other fields. For more details on
variational inequalities, see
\cite{facchinei2003,ferris1997,glowinski1984,glowinski1981,harker1990,kinderlehrer2000,lions1967} and the references therein.

A very popular and successful strategy for computing a solution of a variational
inequality is to reformulate it as a \emph{semismooth rootfinding problem}
\cite{chen2000,deluca1996,hintermueller2002,kummer1988,klatte2002,qi1993,ulbrich2011},
where it is possible to do so. That
is, \eqref{eqn:vi} is equivalently reformulated as the task of finding $z \in Z$
such that
\begin{equation} \label{eqn:ss}
F(z) = 0,
\end{equation}
for a residual function $F: Z \to V$, where $Z$ and $V$ are real Banach spaces
and $Z$ is reflexive. The space $Z$ is typically constructed via $Z = U \times
\Lambda$, where $\Lambda$ is a suitable space of Lagrange multipliers. The
residual $F$ may not be differentiable in the classical Fr\'echet sense but enjoys a
weaker property called semismoothness (defined later in Definition \ref{defn:ss}). The problem \eqref{eqn:ss} is
constructed in such a way that there is a bijection between solutions of
\eqref{eqn:vi} and roots of \eqref{eqn:ss}.
While the
standard Newton--Kantorovich iteration
\begin{equation}
z_{k+1} = z_k - F'(z_k)^{-1} F(z_k)
\end{equation}
requires the existence of the Fr\'echet derivative $F'$ of $F$, it is
possible to define a semismooth Newton iteration for semismooth residuals (defined
later in section \ref{sec:semismooth}). This
method
exhibits locally superlinear convergence under certain regularity conditions on
the solution $z$.

Variational inequalities often admit multiple solutions, and these are typically
significant for the application at hand. For example, a nonconvex optimization
problem may permit several local minima, while a game may permit
multiple Nash equilibria. Identifying these distinct solutions is important for
understanding the system as a whole. The question of calculating distinct
roots of semismooth residuals such as \eqref{eqn:ss} naturally follows. In this
paper, we analyze a numerical
technique called deflation that can successfully identify multiple solutions of
variational inequalities, provided they exist and are isolated from each other.

The central idea of deflation is to compute distinct roots of the semismooth
residual $F$ in the following manner. Let us suppose we are given a semismooth
residual $F: Z \to V$ and a single known root $r \in Z$ that satisfies
some regularity conditions to be made precise later.  Deflation constructs a
modified residual $G: Z \to V$ with the following properties:
\begin{enumerate}
\item[(1)] $G$ preserves roots: $G({z}) = 0 \iff F({z}) = 0$, for ${z} \in
Z\setminus\{r\}$,
\item[(2)] A semismooth Newton method applied to $G$ from any initial
guess in $Z\setminus\{r\}$ will not converge to $r$.
\end{enumerate}
This latter property holds even if semismooth Newton on $G$ is initialized from
the same initial guess that led to the convergence to $r$ in the first instance.
That is, if semismooth Newton is applied to $G$, and it converges, it will
converge to a distinct solution. By enforcing nonconvergence of the
semismooth Newton method to known solutions, deflation enables the discovery of unknown
ones. The deflated problem is constructed via the application of a \emph{deflation
operator} to the underlying problem $F$.

The idea of deflation was first investigated in the context of differentiable maps
$F: \mathbb{R}^n \to \mathbb{R}^n$ by Brown and Gearhart \cite{brown1971}, and
was subsequently reinvented in the context of optimization as the tunneling
method of Levy and G\'omez \cite{levy1985}.
Birkisson \cite{birkisson2013} and Farrell et al.~\cite{farrell2014} analyzed it
in the context of Fr\'echet differentiable maps between Banach spaces, allowing
for its application to smooth partial differential equations. The main
contribution of this work is to extend the theory to the case where $F$ is
semismooth, but not Fr\'echet differentiable.

The importance of multiple solutions of variational inequalities has motivated
other authors to develop various approaches for computing them. A simple
strategy is to vary the initial guess given to the solver \cite{vonstengel2002},
but this is heuristic and labour-intensive \cite{tinloi2003}. By contrast, the
deflation technique does not involve modifying the initial guess; the problem
itself is modified. Judice and Mitra \cite{judice1988} develop an algorithm for
enumerating the solutions of linear complementarity problems, variational inequalities
of the form: find $x \in \mathbb{R}^n_+$ such that
\begin{equation}
(Mx + q)^T (y - x) \ge 0 \quad \text{for all } y \in \mathbb{R}^n_+,
\end{equation}
where $\mathbb{R}^n_+$ is the non-negative orthant of Euclidean space,
$M \in \mathbb{R}^{n \times n}$, and $q \in \mathbb{R}$.
Their algorithm requires exhaustive
exploration of a binary tree whose size is exponential in terms of the size of
the problem, and is thus impractical for large problems. Tin-Loi and Tseng
\cite{tinloi2003} develop an algorithm for finding multiple solutions of linear
complementarity problems by augmenting the problem with constraints that
eliminate known solutions; while very successful on the problems considered, the
size of each problem increases with each solution eliminated. By contrast, the
deflation technique does not increase the size of the problems to be solved
after each solution found. Another strategy is to extend classical
path-following algorithms to parameter-dependent variational inequalities, and
trace out the bifurcation diagram as the parameter is varied
\cite{conrad1988,mittelmann1983}. By continuing around turning points, distinct
solutions for the same parameter values may be identified; however, this
strategy will only identify distinct solutions that happen to lie on a connected
branch. The deflation technique enables the discovery of solutions on
disconnected branches.

Deflation was first applied in the context of finite-dimensional mixed complementarity
problems (a particular kind of variational inequality) by Kanzow
\cite{kanzow2000}.  Kanzow reports some success on rather difficult test
problems, but remarks that deflation is ``not \dots very reliable for larger
problems''. This impression is more widely shared: Allgower \& Georg
\cite{allgower2003} remark (in the context of nonlinear equations, not variational inequalities)
that ``it is often a matter of seeming chance whether one obtains an additional
solution''. We hypothesize that these negative experiences are a consequence of using
a poor deflation operator, the norm deflation operator proposed by Brown \&
Gearhart. We will demonstrate that deflation is much more robust and
effective for semismooth problems with an elementary modification to the deflation operator that
recovers the correct behaviour of the deflated problem at infinity.

This paper is laid out as follows. In section \ref{sec:deflation} we define a
deflation operator, and prove its effectiveness in the semismooth case. The
regularity conditions required on the known solution $z$ are exactly those used
to prove locally superlinear convergence of the semismooth Newton method itself
in \cite{chen2000,hintermueller2002}; no additional assumptions are required. In
section \ref{sec:findimexamples}, we apply the technique to calculating
distinct solutions of several illustrative finite-dimensional variational inequalities, while in section \ref{sec:infdimexamples} we
apply the technique to infinite-dimensional problems with a mesh-independent function-space-based
algorithm.
These
examples demonstrate that the shifted deflation operator applied in this work is
effective in numerical practice. We conclude with some remarks and open
questions in section \ref{sec:conclusion}.

%Continuation methods for calculating a single solution: \cite{kojima1991,chen1993,chen1995,kanzow1996,kanzow1998}

\section{Deflation} \label{sec:deflation}

\subsection{Deflation operators}

Given a known solution $r$, the deflated problem $G(z) = 0$ is constructed by
the application of a \emph{deflation operator} to the original problem $F(z) = 0$:
\begin{equation} \label{eqn:deflatedproblem}
G(z) = M(z; r) F(z).
\end{equation}
The requirements for $M$ to enable the discovery of new solutions are captured in
the definition of a deflation operator. The definition has been slightly modified
from \cite{farrell2014} to allow it to apply to semismooth problems in the sequel.
\begin{definition}[Deflation operator for isolated solutions]
Let $F: D \subset Z \to V$, with its domain of definition $D$ an open subset of $Z$. Let $r \in D$ be an isolated root of $F$, i.e.~$F(r) = 0$ and there exists an open
ball around $r$ with no other roots of $F$.
We say that $M(\cdot, r): Z \setminus \{r\} \to L(V, V)$ is a deflation operator for $F$ at $r$ if
\begin{enumerate}
\item $M(z;r) \in L(V, V)$ is invertible for all $z \ne r$ in a neighborhood of $r$.
\item The deflated residual does not converge to zero as $z \to r$:
\begin{equation} \label{eqn:deflation}
\liminf_{z \to r} \|M(z;r) F(z)\|_{V} > 0.
\end{equation}
\end{enumerate}
\end{definition} Property \eqref{eqn:deflation} is referred to as the deflation property.
The fundamental example of a deflation operator, proposed by Brown \& Gearhart
in 1971 \cite{brown1971}, is
\begin{equation} \label{eqn:normdeflation}
M(z; r) = \frac{\mathcal{I}_V}{\|z - r\|^p_Z},
\end{equation}
where $\mathcal{I}_V$ is the identity map on $V$, and $p \ge 1$. The power $p$ controls the rate of blowup as $z \to r$. This is known to be a deflation operator
in the case where $F$ is continuously Fr\'echet differentiable
\cite{farrell2014}. This operator was one of two considered by Kanzow, and the
operator considered by Allgower \& Georg. However, this operator has a major
drawback: as $\|z-r\| \to \infty$ in any direction, $M(z; r) \to 0$. This
often leads to $G(z) \to 0$ as well, depending on the behaviour of $F(z)$ at
infinity.
Farrell et al.~\cite{farrell2014} suggested a simple modification to the deflation operator to recover
the behaviour $M(z; r) \to 1$ and hence $G(z) \to F(z)$ at infinity: the
addition of a shift. The shifted deflation operator is
\begin{equation} \label{eqn:shifteddeflation}
M(z; r) = \left(\frac{1}{\|z-r\|^p_Z} + 1\right) \mathcal{I}_V.
\end{equation}
This is much more effective in numerical practice than
\eqref{eqn:normdeflation}; the incorrect behaviour at infinity
likely accounts for the unsatisfactory performance reported by
Kanzow and Allgower \& Georg. This will be investigated further
in the examples in section \ref{sec:findimexamples}.

\begin{remark}
After the deflated problem \eqref{eqn:deflatedproblem} is constructed, the Newton--Kantorovich or semismooth Newton algorithms will be applied to it, and therefore the differentiability
of the deflation operator should be established. These deflation operators \eqref{eqn:normdeflation} and \eqref{eqn:shifteddeflation} are differentiable away
from $z = r$ if the norm used on the Banach space $Z$ is differentiable away from zero,
i.e.~if the Banach space $Z$ is Fr\'echet smooth. Note that a reflexive Banach space always
admits an equivalent Fr\'echet smooth norm \cite{fry2002}, and hence this requirement is satisfied
after possibly renorming. We therefore assume this property henceforth.
\end{remark}

\begin{remark}
In practice only an approximation $\tilde{r} \approx r$ is available for
use in the deflation operator. The question then arises of how this
approximation affects the computation of the roots of $G$
(e.g.~Wilkinson \cite[pp.~55]{wilkinson1963} considered this issue in
the context of unshifted deflation for polynomial rootfinding). Since
the shifted deflation operator \eqref{eqn:shifteddeflation} satisfies
$M(z; \tilde{r}) \approx 1$ away from $\tilde{r}$ (and $r$), no
difficulties are encountered for solutions that are sufficiently far
apart. If two solutions are very close together, a simple remedy
discussed by Wilkinson is to calculate a root of $G(z)$, then use that
as initial guess for further Newton iterations on $F(z) = 0$.
\end{remark}

\subsection{The Fr\'echet-differentiable case}
For completeness, we state the result arguing that  \eqref{eqn:normdeflation} and
\eqref{eqn:shifteddeflation} are deflation operators in the Fr\'echet-differentiable case.
Incidentally, this result can be proven analogously to the semismooth case as in Theorem \ref{thm:ssn} below, which provides an alternative proof to the one found in \cite{farrell2014}.

\begin{theorem}[Deflation for Fr\'echet differentiable problems \cite{farrell2014}] \label{thm:frechet}
Let $F: D \to V$ be a continuously Fr\'echet differentiable operator with derivative $F': D \to L(Z, V)$, and let $M$ be given by \eqref{eqn:normdeflation} or \eqref{eqn:shifteddeflation}. Let $r \in D$ be an isolated
solution of $F$, i.e.~satisfy $F(r) = 0$ with $F'(r)$ invertible. Then $M$ is a deflation operator for $F$ at $r$.
\end{theorem}

\begin{remark} \label{rem:othernorms}
It may be more convenient to use another norm $\|\cdot\|_X$ in the deflation
operator, provided $Z \hookrightarrow X$. It is also possible to
use a seminorm $|\cdot|_X$, provided
 $\lim_{z \to r} \|T(z)\|_V / |z - r|_X =
0$, where $T(z)$ is the Taylor remainder associated with 
$F(z) = F(r) + F'(r)(z-r) + T(z)$.
\end{remark}

\subsection{The semismooth case} \label{sec:semismooth}
We now consider the semismooth case.
\begin{definition}[Semismoothness \cite{chen2000,hintermueller2002,mifflin1977}]
\label{defn:ss}
Let $Z$ and $V$ be Banach spaces. Let $F: D \subset Z \to V$, where $D$
is an open subset of $Z$.
$F$ is semismooth at $z \in D$ if it is locally Lipschitz continuous at $z$
%,
%the directional derivative $F'(z; \cdot)$ exists, 
and there exists an open neighbourhood
$N \subset D$ containing $z$ with a \emph{Newton derivative}, i.e.\
a mapping
$H: D \to L(Z, V)$ with the property that
\begin{equation} \label{eqn:littleo}
F(z + h) - F(z) - H(z + h)h = o(h)
\end{equation}
for all $z$ in $N$.
\end{definition}
With this Newton derivative,
the semismooth Newton iteration is given by
\begin{equation}
z_{k+1} = z_k - H(z_k)^{-1} F(z_k).
\end{equation}

We now state the main result of this work. The following theorem is novel.
\begin{theorem}[Deflation for semismooth problems]\label{thm:ssn}
Let $Z$ and $V$ be Banach spaces, and let $F: D \subset Z \to V$.
Let $r$ be a root of $F$. Suppose
$F$ is semismooth at $r$ with Newton derivative $H: D \to L(Z, V)$ in a neighbourhood $N \subset Z$
around $r$.
Further assume that $H(z)$ is invertible for all $z \in N$ and $\{\|H(z)^{-1}\|: z \in N\}$ has a finite upper bound $\Gamma$.
Then the operators \eqref{eqn:normdeflation} and \eqref{eqn:shifteddeflation}
are deflation operators for $p \ge 1$.
\end{theorem}

\begin{remark} These are the same assumptions used to prove the locally superlinear convergence
of the semismooth Newton method in \cite{chen2000,hintermueller2002}.
\end{remark}

\begin{proof}
For brevity, define
\begin{equation} \label{eqn:brevity}
M(z; r) = \left(\frac{1}{\|z-r\|^p_Z} + \sigma\right) \mathcal{I}_V,
\end{equation}
with $\sigma = 0$ corresponding to \eqref{eqn:normdeflation} and $\sigma = 1$
corresponding to \eqref{eqn:shifteddeflation}.
%As in the proof of Theorem \ref{thm:frechet}, define $M$ as in \eqref{eqn:brevity}. 
Invertibility of $M(z; r)$ for all $z \neq r$ is obvious.
Consider $z \in N\setminus\{r\}$. Let $\gamma = \Gamma^{-1}$, and define
\begin{equation}
T(z) = F(z) - F(r) - H(z)(z-r).
\end{equation}
As before, $T(z) = o(z-r)$ from the definition of semismoothness. We then have
\begin{align}
\|M(z; r)F(z)\|_V
%= \frac{\| F(z) \|_V }{\| z - r\|^p_Z}
&\ge
\frac{\| F(z) - F(r) - H(z)(z-r) + H(z)(z-r) \|_V}{\| z - r \|^p_Z} - \sigma \|F(z)\|_V\\
&\ge
\frac{| \| H(z)(z - r) \|_V - \| T(z) \|_V |}{\| z - r \|^p_Z} - \sigma \|F(z)\|_V\\
&\ge
\frac{ \|H(z)(z - r) \|_V - \| T(z) \|_V }{\| z - r \|^p_Z} -\sigma \|F(z)\|_V\\
&\ge \frac{\gamma \|z - r\|_Z - \|T(z)\|_V}{\| z - r \|^p_Z} - \sigma \|F(z)\|_V,
\end{align}
Since $\|F(z)\|_V \to 0$ and $\|T(z)\|_V/\|z-r\|_Z \to 0$ as $z \to r$, we have
\begin{equation}
\liminf_{z \to r} \|M(z; r) F(z)\|_V \ge \liminf_{z \to r} \gamma \|z-r\|_Z^{1-p} > 0.
\end{equation}
as required.
%and the argument follows as in Theorem \ref{thm:frechet}.
\end{proof}

It remains to show that the deflated problem \eqref{eqn:deflatedproblem} is in fact
semismooth. This property is verified in the following result.
\begin{lemma}
Let $F: D \subset Z \to V$ be semismooth at $z \in D, z \neq r$ with Newton derivative $H_F$. Let $M$ be given by 
\begin{equation} \label{eqn:brevity2}
M(z; r) = \left(\frac{1}{\|z-r\|^p_Z} + \sigma\right) \mathcal{I}_V,
\end{equation}
for some $\sigma \ge 0$. Then the product $G(z) = M(z; r) F(z)$ is also semismooth
at $D$ with Newton derivative action
\begin{equation}
H_G(z) h = \left(\sigma + \frac{1}{\|z-r\|^p_Z}\right) H_F(z) h - p \frac{\langle z^*, h\rangle_{Z^*, Z}}{\|z - r\|_Z^{p+1}} F(z),
\end{equation}
where $z^* \in Z^*$ is the derivative of the norm $\|\cdot\|_Z$ at $z-r$, i.e.~satisfies
\begin{equation}
\|z^*\|_{Z^*} \le 1, \quad \langle z^*, z-r \rangle_{Z^*, Z} = \|z-r\|_Z.
\end{equation}
\end{lemma} 
\begin{remark}
If $Z$ is a Hilbert space, then the Riesz representation of $z^*$ is ${(z-r)}/{\|z-r\|}$.
\end{remark}
\begin{proof}
This follows from the well-known calculus rules for semismooth and continuously Fr\'echet differentiable mappings, see e.g., \cite{hintermueller2002,ulbrich2011}
\end{proof}

\begin{remark}
Since the deflated problem is semismooth, the usual sufficient conditions guaranteeing
local superlinear convergence may be applied \cite{hintermueller2002} to the
deflated residual.
\end{remark}

\begin{remark}
Since the deflated problem is also semismooth, any devices developed for globalizing
convergence may be applied, such as line search techniques and continuation
e.g.~\cite{kanzow2000,ralph1994}.
\end{remark}

\begin{remark}\label{rem:rall_rheinboldt}
It would be of significant interest to derive sufficient conditions guaranteeing
the convergence of the same initial guess to two distinct solutions via
deflation. Some initial results in this vein in the smooth case are discussed
in \cite{farrell2015d}, where it is shown that repeated applications of the well-known Rall-Rheinboldt global convergence theorem \cite{rall1974,rheinboldt1978} can assure convergence to two distinct solutions starting from the same initial guess.  A Rall-Rheinboldt-type result, as opposed to Newton-Kantorovich, would be essential to prove convergence to multiple solutions in the context of deflation. This is because the Rall-Rheinboldt theorem places conditions on the radii of convergence of the balls centered at the solutions, rather than guaranteeing the existence of a unique solution in a ball around the initial guess. We briefly investigate the limitations of the classical theory in the framework of semismooth equations below, and in doing so we explain the need for a result of Rall-Rheinboldt-type that is native to the semismooth case.

Many infinite-dimensional semismooth equations of interest share a common structure. In particular, due to low multiplier regularity for bound constrained variational problems, one often resorts to a Moreau-Yosida-type approximation and considers a sequence of (semismooth) equations with residuals taking the form
\[
F_{\gamma}(z) := A(z) + \gamma \Phi(z) - f,
\]
where $A$ is a continuously Fr\'echet differentiable operator, $\Phi$ is a semismooth superposition operator, $\gamma > 0$ is a penalty parameter, and $f$ is constant, cf. \cite{hintermueller2006}. For the sake of argument, assume that $\Omega \subset \mathbb R^n$ is a nonempty, open, and bounded set; $Z = H^1(\Omega)$ the usual Sobolev space of $L^2$-functions with weak derivatives in $L^2$, and $\Phi$ is generated by the function $\phi(x) := \max\{0,x\}$, i.e., $\Phi(z)(x) := \max\{0,z(x)\}$. Since $\phi$ and thus $\Phi$ are nonsmooth, we cannot directly employ the arguments in \cite{farrell2015d}. 

However, by smoothing the $\max$-function, we can obtain a further approximation of the original problem that is regular enough to exploit the Rall-Rheinboldt theory. The remaining question is whether the
convergence guarantees for the smooth problem are stable as $\varepsilon \downarrow 0$. Suppose we replace $\max\{0,x\}$ by
\[
(x)^{\varepsilon}_{+} := \left\{\begin{array}{ll}
x - \frac{\varepsilon}{2}, & x \ge \varepsilon,\\
\frac{x^2}{2\varepsilon}, & x \in (0,\varepsilon),\\
0,& x \le 0.
\end{array}\right.
\]
and define $\Phi_{\varepsilon}(z)(x) := (z(x))^{\varepsilon}_{+}$ and $F_{\gamma,\varepsilon}(z) := A(z) + \gamma \Phi_{\varepsilon}(z) - f$. One of the essential ingredients of the sufficient conditions in the Rall-Rheinboldt theorem are the (local) Lipschitz properties of the derivative $F'_{\gamma,\varepsilon}$ at the distinct solutions. In particular, we require an open neighborhood $E_{i,\varepsilon}$ of each solution $z_{i,\varepsilon}$ along with a constant $\omega_{i,\varepsilon} > 0$ such that
\[
\| F'_{\gamma,\varepsilon}(z_{i,\varepsilon})^{-1}\left(F'_{\gamma,\varepsilon}(u) - F'_{\gamma,\varepsilon}(v)\right)\| \le \omega_{i,\varepsilon} \|u - v\|\quad \forall u,v \in E_{i,\varepsilon}.
\]
Whereas the operator $A$ can be assumed to be unproblematic, one readily derives the estimates 
\[
\aligned
\|  [(\cdot)^{\varepsilon}_+]'(u) -  [(\cdot)^{\varepsilon}_+]'(v) \|_{(H^1)^*} 
&\le 1,\\
| [(\cdot)^{\varepsilon}_+]'(u(x)) -  [(\cdot)^{\varepsilon}_+]'(v(x)) | 
&\le \frac{1}{\varepsilon} | u(x) - v(x) | \text{ a.e. } \Omega
\endaligned
\]
for any two functions $u,v \in H^1(\Omega)$. Although the difference of the smoothed  operators is uniformly bounded in $u,v,$ and $\varepsilon \ge 0$, the pointwise relation, which holds as an equality on the set $\left\{x \in \Omega \left| u(x), v(x) \in
(0,\varepsilon) \right.\right\}$, would indicate that any form of
affine-covariant Lipschitz constant $\omega_{i,\varepsilon}$ would unfavorably depend on $\varepsilon$. As a result, the $\omega_{i,\varepsilon}$-dependent radii associated with the balls of convergence for the distinct roots would converge to zero. Additional assumptions on the structure of $E_{i,\varepsilon}$ that would avoid these issues are unrealistic, e.g., suppose $\Omega \subset \mathbb R^1$ so that $H^1(\Omega) \hookrightarrow C(\overline{\Omega})$ and assume that there exists $\eta > 0$ (independent of $\varepsilon$) such that $z_{i,\varepsilon} \le -\eta < 0$.
\end{remark}

\section{Finite-dimensional examples} \label{sec:findimexamples}
We investigate the effectiveness of the deflation approach by applying it to
various semismooth problems in the literature that exhibit distinct solutions.

\subsection{Complementarity problems}
We consider the nonlinear complementarity problem NCP$(F)$: given $F: \mathbb{R}^n \to \mathbb{R}^n$,
find $z \in \mathbb{R}^n$ such that
\begin{equation} \label{eq:ncp}
z \ge 0, \quad F(z) \ge 0, \quad z \perp F(z).
\end{equation}
This is equivalent to the variational inequality
\begin{equation}
\text{find } z \in K \text{ s.t. } F(z)^T (y - z) \ge 0 \text{ for all } y \in K,
\end{equation}
where
\begin{equation}
K = \{z \in \mathbb{R}^n : z \ge 0\},
\end{equation}
with all inequalities understood componentwise. We apply a standard semismooth reformulation
of the problem using the Fischer--Burmeister NCP function
\begin{align}
\phi_{FB}&: \mathbb{R} \times \mathbb{R} \to \mathbb{R} \nonumber \\
\phi_{FB}(a, b) &= \sqrt{a^2 + b^2} - a - b,
\end{align}
which has the property that $\phi_{FB}(a, b) = 0 \iff a \ge 0, b \ge 0, ab = 0$ \cite{fischer1992}.
The nonlinear complementarity problem \eqref{eq:ncp} is equivalent to finding roots of the
semismooth residual $\Phi: \mathbb{R}^n \to \mathbb{R}^n$ defined by
\begin{align}
\Phi_i(z) = \phi_{FB}(z_i, F_i(z)).
\end{align}
This is then solved with a semismooth Newton method \cite{kanzow2000,kummer1992,qi1993}. As
roots $\{r_1, \dots, r_n\}$ are discovered, semismooth Newton is applied to
\begin{equation} \label{eqn:final}
G(z) = M(z; r_1) M(z; r_2) \cdots M(z; r_n) \Phi(z),
\end{equation}
where $M$ is given by \eqref{eqn:shifteddeflation}. That is, the deflation operators
for each solution are concatenated to deflate all known solutions. Unless noted otherwise,
the parameter choice $p = 2$ was used.

\subsection{Kojima and Shindoh (1986)}
This problem was first proposed by Kojima and Shindoh \cite{kojima1986} and
is an NCP with $F:\mathbb{R}^4\rightarrow\mathbb{R}^4$ given by
\begin{align}
F(z)=\left[\begin{array}{l}
3z_1^2 + 2z_1z_2 + 2z_2^2 + z_3 + 3z_4 - 6\\
2z_1^2 + z_2^2   + z_1 + 10z_3 + 2z_4 - 2\\
3z_1^2 + z_1z_2  + 2z_2^2 +  2z_3 + 9z_4 - 9\\
 z_1^2 + 3z_2^2  + 2z_3 + 3z_4 - 3
\end{array}\right].
\end{align}
It admits two solutions,
\begin{align*}
\begin{array}{lclcl}
{z}^{(1)} = [1, 0, 3, 0]^T, & & \multirow{2}{*}{\text{with residuals}} & & F({z}^1)=[0, 31, 0, 4]^T,\\
{z}^{(2)} = [\sqrt{6}/2, 0, 0, 1/2]^T,& & & & F({z}^2)=[0,2+\sqrt{6}/2, 0, 0]^T.
\end{array}
\end{align*}
This problem was used again by Dirkse and Ferris \cite{dirske1995} as an
example of a problem in which classical Newton solvers struggle to find a solution.
This is because one of the two solutions, ${z}^{(2)}$, has a degenerate third
component, i.e. ${z}^{(2)}_3=F_3({z}^{(2)})=0$, and hence does not satisfy
strict complementarity. Another feature of this problem is that the linear
complementarity problem formed through linearization of the residual $F$ around
zero has no solution, causing difficulties for the Josephy--Newton method there
\cite{josephy1979}.

This is a relatively easy problem to solve and deflation with shifting ($\sigma = 1$) successfully finds both
solutions from many initial guesses. We chose initial guess
$[7/10, \dots, 7/10]^T$. With no line search, semismooth
Newton converged to ${z}^{(1)}$ in 7 iterations; after deflation,
semismooth Newton converged to ${z}^{(2)}$ in 12 iterations. By contrast, without
shifting ($\sigma = 0$) deflation did not identify any additional
solutions.

\subsection{Gould (2001)}
This is a nonconvex quadratic programming problem with linear constraints suggested by N.~I.~M.~Gould in an invited lecture to the 19$^{\mathrm{th}}$ biennial
conference on numerical analysis \cite{gould2002}. It is a quadratic minimization problem
with an indefinite Hessian of the form
\begin{align*}
\min\limits_x f(x)=-2(x_1-1/4)^2+2(x_2-1/2)^2,\hspace{15pt}\text{s.t.}\hspace{15pt}\left\{\begin{array}{r}x_1+x_2\leq 1,\\6x_1+2x_2\leq 3,\\x_1,x_2\geq 0.\end{array}\right.
\end{align*}
The first order Karush--Kuhn--Tucker
optimality conditions yield an NCP with residual
\begin{align}
F(z) = \left[\begin{array}{r}-4(x_1-1/4) + 3\lambda_1+\lambda_2\\4(x_2-1/2)+\lambda_1+\lambda_2\\3-6x_1-2x_2\\1-x_1-x_2\end{array}\right],
\end{align}
where $z=[x,\lambda]$, with $\lambda=[\lambda_1,\lambda_2]$ the vector of the
Lagrange multipliers associated with $F_3(z)\geq0$ and $F_4(z)\geq0$
respectively. Note that in this case it is not necessary to use Lagrange
multipliers to enforce $x\geq0$ as this is implicit in the NCP formulation. The
nonconvexity of the function $f$ makes this problem difficult; it attains two
minima with similar functional values and has a saddle point at $x=[1/4,1/2]^T$. 
The central path to be followed by an interior point method is pathological,
with different paths converging to the different minima.

We directly solve the arising NCP with the semismooth Newton method with
deflation with shifting and without
line search.
The initial guess was
$[2/10, 2/10, 0, 0]^T$.
In order, the three solutions found were
\begin{align*}
\begin{array}{lcl}
{z}^{(1)} = [1/4, 1/2, 0, 0]^T,& & F({z}^{(1)})=[0, 0, 1/4, 1/4]^T,\\
{z}^{(2)} = [0, 1/2, 0,0]^T, & \multirow{1}{*}{\text{with residuals}} & F({z}^{(2)})=[1, 0, 1, 1/2]^T,\\
{z}^{(3)} = [11/32, 15/32, 1/8, 0]^T, & & F({z}^{(3)})=[0, 0, 0, 3/16]^T.
\end{array}
\end{align*}
These are the saddle point, the global minimum and the local minimum respectively.
The KKT conditions make no distinction between minima and saddle points, and
hence the solver finds both kinds of stationary points. The number of iterations
required was 5, 7 and 10 respectively. As before, without shifting deflation did not
successfully identify any additional solutions.

\subsection{Aggarwal (1973)}
This is a Nash bimatrix equilibrium problem arising in game theory.
This kind of problem was first introduced by von Neumann and
Morgenstern \cite{vonneumann1945} and the existence of its solutions was
further studied by Nash \cite{nash1951} and Lemke and Howson
\cite{lemke1964}. In the same paper, Lemke and Howson also presented a numerical
algorithm for computing solutions to these kinds of problems. This
example was introduced by Aggarwal \cite{aggarwal1973} to prove that it is
impossible to find all solutions of such problems using a modification of
the Lemke--Howson method that had been conjectured to compute all solutions.

The problem consists of finding the equilibrium points of a bimatrix (non-zero
sum, two person) game. Let $A$ and $B$ be the $n \times n$ payoff matrices of
players $1$ and $2$ respectively. Let us assume that player $1$ plays the
$i^{th}$ pure strategy and player $2$ selects the $j^{th}$ pure strategy amongst
the $n$ strategies available to each. The entries of $A$ and $B$, $a_{i,j}$ and
$b_{i,j}$ respectively, correspond to the payoff received by each player. It is then possible
to define a mixed strategy for a player which consists of a $n \times 1$ vector
$x$ such that $x_i\geq0$ and $x_1+...+x_n=1$. Denote by $x$ and $y$ the mixed
strategies for player $1$ and $2$ respectively. The entries of these vectors
stand for the probability of the player adopting the corresponding pure
strategy. The expected payoffs of the two players are then $x^TAy$ and $x^TBy$
respectively. An equilibrium point $(x^*,y^*)$ is reached when, for all $x$,
$y$,
\begin{align}
(x^*)^TAy^*\geq x^TAy^*,\hspace{12pt}\text{and}\hspace{12pt}(x^*)^TBy^*\geq (x^*)^TBy,
\end{align}
i.e. neither player can unilaterally improve their payoff.

Aggarwal's counterexample admits three Nash equilibria.
These equilibria are related to the solutions of the NCP with residual
\begin{align}
F(z) =
\begin{pmatrix}
\overline{A}y - e \\
\overline{B}^Tx - e
\end{pmatrix}
,
\end{align}
where $z = [x, y]^T$ and $e = [1, 1, \dots, 1]^T$,
$\overline{A}$ and $\overline{B}$ are positive-valued loss matrices related to $A$ and $B$
respectively, and $x$ and $y$ relate to the mixed strategy adopted by each player \cite[\S 1.4]{murty1988}.
The data for this problem is
\begin{align*}
\overline{A} =
\begin{bmatrix}
30 & 20 \\
10 & 25
\end{bmatrix},
\text{ and }
\overline{B} =
\begin{bmatrix}
30 & 10 \\
20 & 25
\end{bmatrix}.
\end{align*}
This problem is quite difficult, and we therefore turned
to continuation to aid convergence, as described below. The problem was modified
to introduce an artificial parameter $\mu$
\begin{align}
F_{\mu}(z) =
\begin{pmatrix}
\mu \overline{A}y - e \\
\mu \overline{B}^Tx - e
\end{pmatrix}
,
\end{align}
with the original problem given by $\mu = 1$. With deflation with
shifting, three solutions
were found for $\mu = 1/1000$ from the initial guess $[0, \dots, 0]^T$, in
5, 24 and 26 iterations of semismooth Newton respectively. (As in the
previous examples, deflation without shifting did not identify any additional
solutions.) All three
branches were then successfully continued to $\mu = 1$ using 50 equispaced
continuation steps and simple zero-order continuation, i.e.~the solution for
the previous value $\mu_-$ is used as initial guess for the solution of the next
value $\mu_+$.
The
three solutions found were
\begin{align*}
\begin{array}{lclcl}
{z}^{(1)} = [0, 1/20, 1/10, 0]^T, & & \multirow{3}{*}{\text{with residuals}} & & F({z}^{(1)})=[2, 0, 0, 1/4]^T,\\
{z}^{(2)} = [1/110, 4/110, 1/110, 4/110]^T,& & & & F({z}^{(2)})=[0, 0, 0, 0]^T,\\
{z}^{(3)} = [1/10, 0, 0, 1/20]^T, & & & & F({z}^{(3)})=[0, 1/4, 2, 0]^T.\\
\end{array}
\end{align*}
Aggarwal observed that the conjectured scheme mentioned above could compute
${z}^{(1)}$ and ${z}^{(3)}$, but could not compute ${z}^{(2)}$.

\subsection{G\'erard, Lecl\`ere and Philpott (2017)} \label{sec:gerard}

G\'erard et al.~describe a stochastic market where the agents are
risk-averse, i.e.~estimate their welfare using a coherent risk measure \cite{gerard2017}.
They give an example of an incomplete market with three different
equilibria, two stable and one unstable. The authors examine
the convergence of the well-known PATH solver \cite{dirske1995b,ferris2000}, and discover that PATH
always yields the unstable equilibrium, even when initialized from
many distinct initial guesses. (An alternative t\^{a}tonnement algorithm does discover
all three equilibria when initialized from different initial guesses.) We therefore investigate
whether deflation can assist a semismooth Newton method in discovering all three
solutions from a single initial guess.

Mathematically, the problem is a mixed complementarity problem, a generalization of
nonlinear complementarity problems. Let $\mathbb{R}_\infty := \mathbb{R} \cup \{-\infty, +\infty\}$. Given $F: \mathbb{R}^N \to \mathbb{R}^N$, a lower bound $l \in \mathbb{R}_\infty^N$ and an upper
bound $u \in \mathbb{R}_\infty^N$, the task is to find $z \in \mathbb{R}^N$ such that exactly
one of the following holds for each $i = 1, \dots, N$:
\begin{enumerate}
\item[(a)] $l_i \le z_i \le u_i$ and $F_i(z) = 0$;
\item[(b)] $l_i = z_i$ and $F_i(z) > 0$;
\item[(c)] $z_i = u_i$ and $F_i(z) < 0$.
\end{enumerate}
This is referred to as MCP($F, l, u$). NCP($F$) is a special case with the particular
choice $l = [0, \dots, 0]^T$ and $u = [\infty, \dots, \infty]^T$. 

The problem at hand is given by $F: \mathbb{R}^{10} \to \mathbb{R}^{10}$, where
\begin{equation}
F(z) = \begin{pmatrix}
          -(\frac{3}{4}(\pi_1 - \frac{23}{2}x_0) + \frac{1}{4}(\pi_2 - \frac{23}{2}x_0))u_4 - (\frac{1}{4}(\pi_1 - \frac{23}{2}x_0) + \frac{3}{4}(\pi_2 - \frac{23}{2}x_0))u_5 \\
          (-\frac{3}{4}(\pi_1 - x_{11}))u_4 + (-\frac{1}{4}(\pi_1 - x_{11}))u_5 \\ 
          (-\frac{1}{4}(\pi_2 - \frac{7}{2}x_{12}))u_4 + (-\frac{3}{4}(\pi_2 - \frac{7}{2}x_{12}))u_5 \\ 
          -(4 - \pi_1 - 2y_1) \\ 
          -(9.6 - \pi_2 - 10y_2) \\ 
          x_0 + x_{11} - y_1 \\ 
          x_0 + x_{12} - y_2 \\ 
          \frac{3}{4}(\pi_1(x_0 + x_{11}) - \frac{23}{4}x_0^2 - \frac{1}{2}x_{11}^2) + \frac{1}{4}(\pi_2(x_0 + x_{12}) - \frac{23}{4}x_0^2 - \frac{7}{4}x_{12}^2) - \theta_P \\ 
          \frac{1}{4}(\pi_1(x_0 + x_{11}) - \frac{23}{4}x_0^2 - \frac{1}{2}x_{11}^2) + \frac{3}{4}(\pi_2(x_0 + x_{12}) - \frac{23}{4}x_0^2 - \frac{7}{4}x_{12}^2) - \theta_P \\
          u_4 + u_5 - 1
\end{pmatrix}.
\end{equation}
with $z = (x_0, x_{11}, x_{12}, y_1, y_2, \pi_1, \pi_2, u_4, u_5, \theta_P)$. The
bounds are given by $l = [0, \dots 0, -\infty]^T$ and $u = [\infty, \dots, \infty]^T$,
i.e.~all variables except $\theta_P$ have lower bound 0, and $\theta_P$ is unconstrained.

To demonstrate that the deflation concept is not confined to a particular
semismooth reformulation, in this example we use an alternative
NCP function. Define
\begin{align} \label{eqn:phiMP}
\phi_{\text{MP}}&: \mathbb{R} \times \mathbb{R} \to \mathbb{R} \nonumber \\
\phi_{\text{MP}}(a, b) &= b - \max{(0, b - a)}
\end{align}
which again has the property that $\phi_{\text{MP}}(a, b) = 0 \iff a \ge 0, b \ge 0, ab = 0$.
The semismooth reformulation of the MCP employed is
\begin{equation}
\Phi_i(z) = \begin{cases}
\phi_{\text{MP}}(z_i, F_i(z)) & i = 1, \dots 9 \\
F_i(z)                 & i = 10.
\end{cases}
\end{equation}

Deflation with shifting was applied from the initial guess $z_0 = [0, \dots, 0]^T$ with $p = 1$,
and a line search algorithm was used to aid convergence (Alg.~2 of Brune et al.~\cite{brune2015}).
With these parameters, the procedure identified all three solutions. With $p = 2$ or
without the line search, only two solutions were found. In order, the three
solutions found were
\begin{flalign}
& [{\pi_1}, {\pi_2}]^{(1)} = [1.2256, 2.0698],\\
& [{\pi_1}, {\pi_2}]^{(2)} = [1.2478, 2.1564],\\
& [{\pi_1}, {\pi_2}]^{(3)} = [1.2358, 2.1095],
\end{flalign}
where only the equilibrium prices are shown for brevity. The solutions were found with 15, 9 and 17
semismooth Newton iterations respectively. The solution found by PATH is the latter. As
with the previous examples, deflation without shifting did not identify any additional
solutions.

\section{Infinite-dimensional examples} \label{sec:infdimexamples}
\subsection{Solving infinite-dimensional variational inequalities}
When solving inequality-constrained infinite-dimensional problems, additional care must be
taken. The main issue here is a general lack of regularity of the Lagrange
multipliers. As a result, it is often impossible to derive a complementarity
system (analogous to KKT-conditions for nonlinear programs) that can be
reformulated as a single semismooth equation. Even in situations where the
associated multiplier is regular enough to allow such a reformulation, we
encounter insurmountable issues in the derivation of a function-space-based
generalized Newton method.  Ignoring these issues and taking a
first-discretize-then-optimize approach will generally lead to mesh-dependent
convergence, i.e.~the number of iterations required to converge increases
significantly as the mesh is refined.

This is illustrated in the following example. Suppose $\Omega \subset \mathbb
R^n$ is a nonempty, open, and bounded subset and let $\mathcal{J}\!:\!H^1_0(\Omega) \to \mathbb R$ be G\^{a}teaux differentiable. Consider the model
problem
\begin{equation}\label{eq:nonconvex}
\min \left\{ \mathcal{J}(u) \text{ over } u \in H^1_0(\Omega) \left|\; u(x) \ge 0 \text{ a.e. } x \in \Omega \right.\right\}.
\end{equation}
We denote the feasible set by $K$.  If \eqref{eq:nonconvex} admits a solution
$\overline{u}$, then we have
\[
\mathcal{J}'(\overline{u})(v - \overline{u}) \ge 0 \quad \forall v \in K
\]
with $\mathcal{J}'(\overline{u}) \in H^{-1}(\Omega)$.
Since $K$ is a cone, an equivalent formulation holds:
\begin{equation}\label{eq:comp_prob}
\mathcal{J}'(\bar{u}) + \lambda = 0,\quad
\overline{u} \in K,\quad \lambda \in K^{\circ},\quad \langle \lambda,\overline{u} \rangle = 0,
\end{equation}
where $K^{\circ}$ is the polar cone to $K$ given by
\[
K^{\circ} := \left\{ v \in H^{-1}(\Omega) \left|\; \langle v ,\varphi \rangle \le 0\; \forall \varphi \in H^1(\Omega) : \varphi \ge 0\right.\right\}.
\]
According to the Radon--Riesz theorem, $\lambda \in K^{\circ}$ is in fact a locally
finite Radon measure on $\Omega$ \cite[pg.~564]{bonnans2000}. Therefore, $\lambda$
cannot in general be evaluated pointwise, in which case \eqref{eq:comp_prob} cannot be reformulated
as a semismooth system of equations.

Suppose further that $\Omega$ is a convex polyhedron
and that $\mathcal{J}'({u}) = A{u} - f$, with $A$ a second-order
linear elliptic operator with smooth coefficients and $f \in L^2(\Omega)$.
Then $\overline{u} \in H^2(\Omega) \cap H^1_0(\Omega)$
\cite[chap.~IV]{kinderlehrer2000} and thus $\lambda \in L^2(\Omega)$.
We may then rewrite \eqref{eq:comp_prob} as
 \begin{equation}\label{eq:infncp}
\mathcal{J}'(\bar{u}) + \lambda = 0,\quad
\lambda = (\lambda - \overline{u})_+,
\end{equation}
where $(x)_+ := \max\{0, x\}$.
Even in this ideal case, the
nonsmooth superposition operator $\Phi(\lambda,u) := \lambda - (\lambda - {u})_+$ must be defined from
$L^2(\Omega) \times H^1_0(\Omega$) into $L^2(\Omega)$. The natural choice for
the generalized derivative of $\Phi$ is given by
\[
\mathcal{G}(\lambda,{u})(\delta \lambda, \delta u) = \chi_{\{\lambda - {u} > 0 \}}\delta u + \chi_{\{ \lambda - {u} \le 0\}} \delta \lambda.
\]
In order for this to be a Newton derivative, the following approximation property must hold:
\[
\|\Phi(\lambda + \delta \lambda, u + \delta u) - \Phi(\lambda, u) - \mathcal{G}(\lambda + \delta \lambda, u + \delta u)(\delta \lambda, \delta u)\| = o(\|(\delta \lambda, \delta u)\|).
\]
However, this only holds true if $\Phi$ is defined from $L^{2+\varepsilon}(\Omega) \times H^1_0(\Omega) \to L^2(\Omega)$
for $\varepsilon > 0$ \cite{hintermueller2002,ulbrich2011}, which is not available even in this ideal case. This so-called ``missing norm gap''
persists for all other known NCP-functions. As a result, the infinite-dimensional problem \eqref{eq:infncp} is not semismooth, which manifests
itself as mesh-dependence on the discrete level \cite{hintermueller2004}.

An alternative mesh-independent scheme can be constructed from the Moreau--Yosida
regularization of the
indicator functional for the constraints with respect to the $L^2(\Omega)$ topology. The key property of this scheme is that the plus
function $(\cdot)_+$ is only applied to $u$, and since $u \in H^1_0(\Omega) \hookrightarrow L^{2+\varepsilon}(\Omega)$,
a norm gap holds and this operator is semismooth.

We sketch the approach taken in our implementation; for more details, see \cite{hintermueller2006,MHintermueller_KKunisch_2006}. Let
\begin{equation}
\mathcal{I}(u) := \int_\Omega i_{\mathbb{R}_+}\! \left(u(x)\right)   \mathrm{d}x,
\end{equation}
where $i_{\mathbb{R}_+}(x)$ is the indicator function for $\mathbb{R}_+$ (i.e.~0 if $x \le 0$ and $+\infty$ otherwise).
The minimization problem \eqref{eq:nonconvex} is equivalent to the unconstrained problem
\begin{equation}\label{eq:nonconvexindicator}
\min \left\{ \mathcal{J}(u) + \mathcal{I}(u) \text{ over } u \in H^1_0(\Omega) \right\}.
\end{equation}
Approximating $\mathcal{I}(u)$ by its $L^2(\Omega)$ Moreau--Yosida regularization
\begin{equation}
\mathcal{I}_\gamma(u) := \inf_{v \in L^2(\Omega)} \{ \mathcal{I}(v) + \frac{\gamma}{2} \|u-v\|^2_{L^2(\Omega)}\} = \frac{\gamma}{2} \int_\Omega (-u)_+^2 \ \mathrm{d}x,
\end{equation}
we obtain a sequence of $\gamma$-dependent problems of the form
\begin{equation}\label{eq:mysubproblem}
\min \left\{ \mathcal{J}(u) + \mathcal{I}_\gamma(u) \text{ over } u \in H^1_0(\Omega) \right\},
\end{equation}
with penalty parameter $\gamma \to \infty$.
The associated first-order necessary condition is
 \begin{equation}\label{eq:ncp_gamma}
\mathcal{J}'(u) - \gamma(-u)_+ = 0,
\end{equation}
which is semismooth for the reasons outlined above. An initial $\gamma$ is
chosen and $\overline{u}_{\gamma}$ computed.  Once this is found, the solver
continues in $\gamma$, using an analytical path-following scheme to drive the
penalty parameter $\gamma \to \infty$ efficiently
\cite{hintermueller2006,MHintermueller_KKunisch_2006,LAdam_HHintermueller_TMSurowiec_2018}.
The mesh and $\gamma$ are linked; as $\gamma \to \infty$, the mesh is uniformly
refined to ensure balanced error estimates (cf.~\cite{hintermueller2009}). At every update of $\gamma$, the
mesh is refined zero or more times until
\begin{equation}
h \le \frac{1}{\sqrt{\gamma}}
\end{equation}
is satisfied, where $h$ is the characteristic mesh size. The process terminates once $\gamma$ reaches a target
value $\gamma_{\max}$, which in this work is taken to be $\gamma_{\max} = 10^6$.

\subsection{Computing multiple solutions of infinite-dimensional variational inequalities with deflation} \label{sec:infdimdefl}
Deflation can be combined with the Moreau--Yosida solver as follows. For the
initial value of $\gamma$, a set of initial guesses is supplied. Each guess is
used as a starting point for semismooth Newton applied to \eqref{eq:ncp_gamma};
if the guess is successful, the solution found is deflated, and the guess is
attempted again. When all guesses have been exhausted, the analytical
path-following strategy is applied to all solutions found, and the next value of
$\gamma$ is taken to be the minimum of these (the most conservative update of
all solutions found). If necessary, the mesh is refined and all solutions
prolonged. The solutions found for the previous step are then used as initial
guesses for the next, until the process terminates with $\gamma \ge \gamma_{\max}$.

\subsection{Zeidler (1988)} \label{sec:zeidler}
Zeidler \cite[pg.~320]{zeidler1988} and Ph\'u \cite{phu1987} study a long thin elastic rod under the action
of a compressive load constrained to lie in a channel of fixed width. Let $s \in [0, L]$ denote the arclength of the rod,
$y \in H^1_0(0, L)$ denote its vertical displacement from centerline of the channel, and $\theta \in H^1(0, L)$ denote the angle between the rod
and the centerline of the channel. The potential energy of the system
is given by
\begin{equation}
\mathcal{E}(y, \theta) = \int_0^L B \left(\theta'(s)\right)^2 + P \cos{\theta(s)} - P - \rho g y \ \mathrm{d}s,
\end{equation}
where $B \in \mathbb{R}$ is the bending stiffness of the rod, $P \in \mathbb{R}$ is the compressive load applied to the right end-point,
$\rho \in \mathbb{R}$ is the mass per unit length of the rod, and $g \in \mathbb{R}$ is the acceleration due to gravity.
The rod is placed in a channel of width
$2\alpha$ such that
\begin{equation}
y(s) \in [-\alpha, \alpha]
\end{equation}
is satisfied almost everywhere. Equilibria of the system are therefore given
by the local minimizers of
\begin{equation} \label{eqn:nonlinearenergy}
\begin{aligned}
& \underset{y \in H^1_0(0, L),\ \theta \in H^1(0, L)}{\text{minimize}}
& & \mathcal{E}(y, \theta) \\
& \underset{\smash{\phantom{y \in H^1_0(0, L),\ \theta \in H^1(0, L)}}}{\text{subject to}}
& & \sin(\theta) = y',\\
& & & \left|y\right| \le \alpha \ \text{a.e.}
\end{aligned}
\end{equation}
Zeidler and Ph\'u consider small deformations and linearize the system using the Taylor expansions
\begin{equation} \label{eqn:zeidlertaylor}
\sin(\theta) \approx \theta, \quad 
\cos(\theta) \approx 1 - \frac{\theta^2}{2}
\end{equation}
yielding the system
\begin{equation} \label{eqn:linearized}
\begin{aligned}
& \underset{y \in H^1_0(0, L),\ \theta \in H^1(0, L)}{\text{minimize}}
& & E(y, \theta) = \int_0^L B (\theta')^2 - P \theta^2 - \rho g y \ \mathrm{d}s \\
& \underset{\smash{\phantom{y \in H^1_0(0, L),\ \theta \in H^1(0, L)}}}{\text{subject to}}
& & \theta = y',\\
& & & \left|y\right| \le \alpha \ \text{a.e.}
\end{aligned}
\end{equation}
where irrelevant constant terms in the functional have been neglected.
Substituting the constraint $\theta = y'$, we arrive at the final system
\begin{equation} \label{eqn:linearized_H2}
\begin{aligned}
& \underset{y \in H^2(0, L) \cap H^1_0(0, L)}{\text{minimize}}
& & J(y) = \int_0^L B (y'')^2 - P (y')^2 - \rho g y \ \mathrm{d}s \\
& \underset{\smash{\phantom{y \in H^2(0, L) \cap H^1_0(0, L)}}}{\text{subject to}}
& & \left|y\right| \le \alpha \ \text{a.e.}
\end{aligned}
\end{equation}

In the absence of the inequality constraint on $y$, the optimality conditions
for this linearized problem comprise a linear beam equation, with either a
unique solution or a one-dimensional nullspace. However, in the presence of the
inequality constraint $|y| \le \alpha$, the problem remains nonlinear and can
support distinct isolated solutions.

The system \eqref{eqn:linearized_H2} was discretized using $H^2(0,
L)$-conforming cubic Hermite finite elements \cite{kirby2018a,kirby2018b} using
Firedrake \cite{rathgeber2016}. The arising linear systems were solved using
the sparse LU factorization of PETSc \cite{balay2017}.

In the absence of gravity ($g = 0$), Ph\'u proved that the first bifurcation
of the system \eqref{eqn:linearized} occurs at $P = B \pi^2/L^2$. In the
presence of small gravity, the reflective symmetry of the system is broken and the zero
solution is no longer a trivial solution, but the bifurcation point will
be nearby. We therefore consider the system
for parameter values $B = 1, g = 1, \rho = 1, L = 1, \alpha = 0.4$ and $P = 10.4$. This
choice of $P$ is sufficiently greater than $B \pi^2/L^2 \approx 9.87$ that it
is reasonable to expect the system to support distinct solutions.

\begin{figure}
\centering
\subcaptionbox{\label{fig:zeidlera}}{
\scalebox{0.3}{\includegraphics[width=\textwidth]{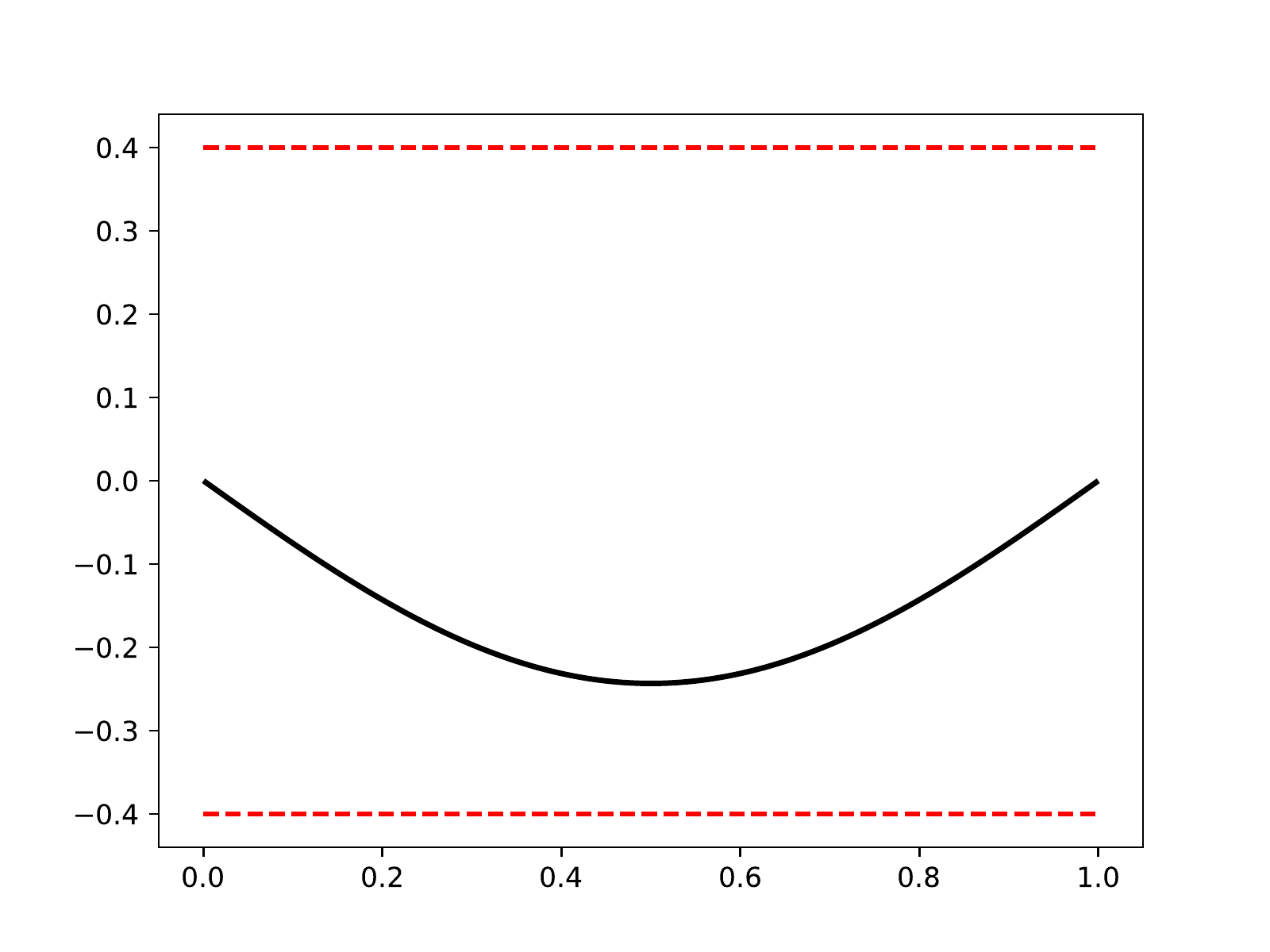}
}}
\
\subcaptionbox{\label{fig:zeidlerb}}{
\scalebox{0.3}{\includegraphics[width=\textwidth]{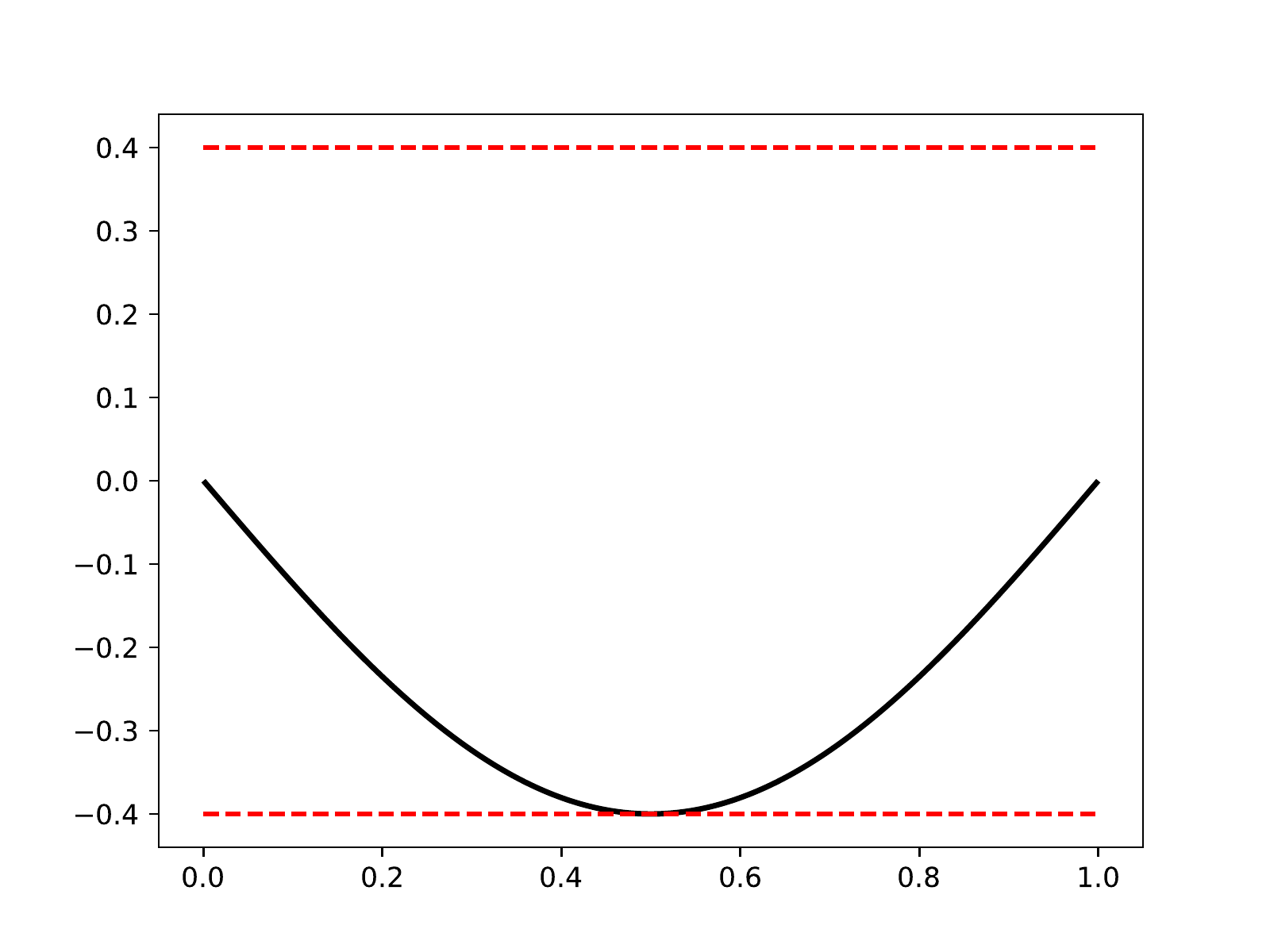}
}}
\
\subcaptionbox{\label{fig:zeidlerc}}{
\scalebox{0.3}{\includegraphics[width=\textwidth]{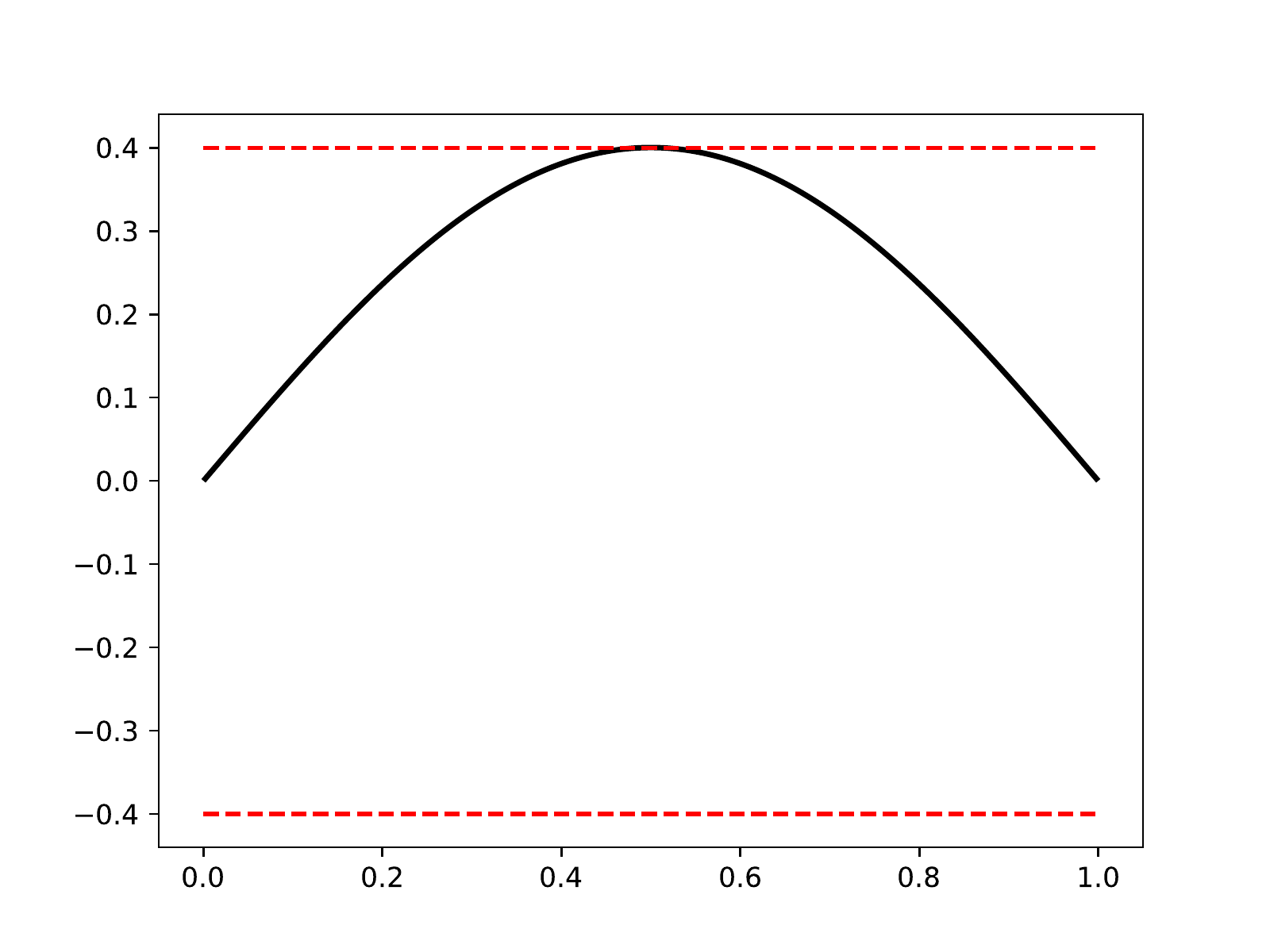}
}}
\caption{Solutions of the linearized beam problem \eqref{eqn:linearized_H2}. The
dashed red lines denote the inequality constraints on the vertical displacement
of the beam.}
\end{figure}

The mesh-independent Moreau--Yosida solver with deflation was applied
with initial guess $y = 0$ for $\gamma = 10$.
This converged in one iteration to the first solution, Figure
\ref{fig:zeidlera}. This is expected as the inequality constraints are inactive
at this solution and the problem is therefore equivalent to the linear beam
equation. Deflation with shifting was then applied with deflation operator
\begin{equation} \label{eqn:zeidlerdefop}
M\left( y, r_y \right) = \|y - r_y\|^{-2}_{L^2(0, L)} + 1,
\end{equation}
where $r_y$ denotes the solution already known. The solver was reinitialized from
the zero initial guess and converged after 6 semismooth Newton iterations to
a second solution that violates the lower bound. (For this low value of $\gamma$, the
bound constraints are only weakly enforced.)
This solution was then deflated using the same operator \eqref{eqn:zeidlerdefop}
and the solver was re-initialized from the zero initial guess. The procedure then converged
after 14 semismooth Newton iterations to the third solution that violates the upper bound\footnote{We also experimented with a mesh-dependent semismooth Newton method applied to this problem. The second solution was found after 42 iterations, while the third was found after 45, and both required a line search. In this case the convergence of the mesh-independent scheme is much more robust.}.
These three solutions were then continued to $\gamma = \gamma_{\max}$ in 9 continuation
steps, with no further solutions found. The three solutions found for $\gamma = \gamma_{\max}$
are shown in Figures \ref{fig:zeidlera}--\ref{fig:zeidlerc}.

This experiment demonstrates an important property of the deflation strategy:
deflation is capable of computing distinct solutions of infinite-dimensional
variational inequalities whose solutions exhibit both nontrivial active sets and
no activity whatsoever, from the same initial guess.

\subsection{A two-dimensional beam under axial compression with obstacle constraints} \label{sec:hyper2d}
In this example we consider a two-dimensional analogue of the previous problem,
and compute several equilibrium configurations of a hyperelastic beam under axial
compression with obstacle constraints.

The physical model employed is compressible neo-Hookean hyperelasticity. Let
$\Omega = (0, 1) \times (0, 1/10)$ denote the undeformed reference configuration,
with boundary $\partial \Omega = \partial \Omega_{\textrm{left}} \cup \partial
\Omega_{\textrm{bottom}} \cup \partial \Omega_{\textrm{right}} \cup \partial
\Omega_{\textrm{top}}$.
Homogeneous Dirichlet conditions are imposed on $\partial \Omega_{\textrm{left}}$,
axial compression Dirichlet conditions are imposed on $\partial \Omega_{\textrm{right}}$,
and natural boundary conditions are imposed on $\partial \Omega_{\textrm{top}}$ and $\partial \Omega_{\textrm{bottom}}$.

In addition, box constraints are imposed on the vertical component $u_2$ of the displacement vector
field $u: \Omega \to \mathbb{R}^2$
on $\partial \Omega_{\textrm{top}}$ and $\partial \Omega_{\textrm{bottom}}$. Let
\begin{equation}
\tau_{\mathrm{top}}(u) = \mathrm{tr}_{\partial \Omega_{\mathrm{top}}} (u_2)
\end{equation}
where $\mathrm{tr}_{\Gamma}: H^1(\Omega) \to H^{1/2}(\Gamma)$ is
the standard trace operator, and let $\tau_{\mathrm{bottom}}$ be defined analogously.
Since $H^{1/2}(\Gamma) \hookrightarrow L^2(\Gamma)$,
we may impose pointwise bound constraints of the type $\tau_{\mathrm{top}}(u) \le \alpha$ and $\tau_{\mathrm{bottom}}(u) \ge \alpha$.

For a given axial compression $\varepsilon$, we seek a displacement vector
\begin{equation}
u \in V_\varepsilon := 
\{v \in H^1(\Omega; \mathbb{R}^2) : \left.v\right|_{\partial \Omega_{\textrm{left}}} = (0, 0), \left.v\right|_{\partial \Omega_{\textrm{right}}} = (-\varepsilon, 0)\}
\end{equation}
that satisfies
\begin{equation} \label{eqn:hyper2d}
\begin{aligned}
&\underset{u \in V_\varepsilon}{\text{minimize}}\;
J(u) = \int_\Omega \psi(u) \ \mathrm{d}x - \int_\Omega B\cdot u\ \mathrm{d}x \\
&\begin{array}{lrcl}
{\text{subject to}}& \tau_{\mathrm{top}}(u)           &\le&     \phantom{-}\alpha \ \text{a.e.} \\
              & \tau_{\mathrm{bottom}}(u) &\ge&     -\alpha \ \text{a.e.},
\end{array}
\end{aligned}
\end{equation}
where $\psi(u)$ is the isotropic compressible neo--Hookean strain energy density with
Young's modulus $E = 10^6$ and Poisson ratio $\nu = 0.3$, $B = (0, -1000)$ is the body force
density due to gravity, and $\alpha$ is the value of the bound constraint
enforced. The problem is discretized using piecewise linear finite elements;
the coarsest grid employed has 3200 triangular elements.

The goal is to solve this problem for $\varepsilon = 0.15$ with $\alpha = 8
\times 10^{-2}$. To do this, continuation is employed. The Moreau--Yosida
regularization of \eqref{eqn:hyper2d} is solved with fixed $\gamma = 100$ from
$\varepsilon = 10^{-3}$ to $\varepsilon = 0.15$ with steps of $\Delta
\varepsilon = 10^{-3}$, with deflation employed at each continuation step. This
initial continuation process yields seven solutions at $\varepsilon = 0.15$,
three inactive solutions and four active solutions. The continuation procedure
with mesh refinement described in section \ref{sec:infdimdefl} is then applied
to continue these solutions from $\gamma = 100$ to $\gamma = 10^6$. On the
coarsest grid, the LU algorithm of MUMPS \cite{amestoy2001} is used to solve the
linear systems arising in semismooth Newton; on finer meshes GMRES-accelerated
geometric multigrid is employed, using all levels in the hierarchy, with three iterations of
Chebyshev-accelerated point-block SOR as a smoother (see Ulbrich et al.~\cite{ulbrich2017}
for rigorous analysis of multigrid in a Moreau--Yosida regularization context). For all solves, the full
undamped semismooth Newton step is used, i.e.~no line search is found to be
necessary for this problem.

\begin{table}
\centering
\begin{tabular}{ccc}
\toprule
Solution & Discovered at & Semismooth Newton iterations \\
\midrule
1 & $\varepsilon = 0$ & 3 \\
2 & $\varepsilon = 4.0 \times 10^{-2}$ & 13 \\
3 & $\varepsilon = 4.1 \times 10^{-2}$ & 9 \\
4 & $\varepsilon = 7.0 \times 10^{-2}$ & 7 \\
5 & $\varepsilon = 7.1 \times 10^{-2}$ & 13 \\
6 & $\varepsilon = 1.44 \times 10^{-1}$ & 16 \\
7 & $\varepsilon = 1.45 \times 10^{-1}$ & 17 \\
\bottomrule
\end{tabular}
\caption{The values of $\varepsilon$ at which new solutions were discovered via
deflation, along with the number of semismooth Newton iterations required.}
\label{tab:discoveryitercounts}
\end{table}

\begin{table}
\centering
\begin{tabular}{ccccc}
\toprule
$\gamma$ & \small{\# Refs.} & \small{Dofs} & \small{Avg SSN its} & \small{Avg GMRES/MG its per SSN step}\\
\midrule
$1.00 \times 10^2$ & 0 & $3.36 \times 10^3$ & 3.28 & -\\
$1.33 \times 10^2$ & 0 & $3.36 \times 10^3$ & 4.00 & -\\
$2.23 \times 10^2$ & 0 & $3.36 \times 10^3$ & 4.50 & -\\
$4.47 \times 10^2$ & 0 & $3.36 \times 10^3$ & 5.00 & -\\
$1.05 \times 10^3$ & 0 & $3.36 \times 10^3$ & 4.25 & -\\
$2.80 \times 10^3$ & 1 & $1.31 \times 10^4$ & 4.14 & 19.03 \\
$8.42 \times 10^3$ & 2 & $5.18 \times 10^4$ & 3.85 & 19.63 \\
$2.81 \times 10^4$ & 3 & $2.06 \times 10^5$ & 3.85 & 22.30 \\
$1.03 \times 10^5$ & 4 & $8.21 \times 10^5$ & 4.14 & 24.62\\
$4.13 \times 10^5$ & 5 & $3.28 \times 10^6$ & 5.00 & 26.17\\
$1.00 \times 10^6$ & 5 & $3.28 \times 10^6$ & 4.00 & 25.75\\
\bottomrule
\end{tabular}
\caption{Number of mesh refinements, degrees of freedom, nonlinear and linear iteration counts required for the continuation in $\gamma$ along $\varepsilon = 0.15$. The
semismooth Newton solver exhibits $\gamma$- and mesh-independence, while the number of Krylov iterations
per semismooth Newton step grows very slowly as $\gamma$ and $h$ are refined.}
\label{tab:krylovitercounts}
\end{table}

\begin{figure}
\begin{tabular}{c}
\includegraphics[width=12.0cm]{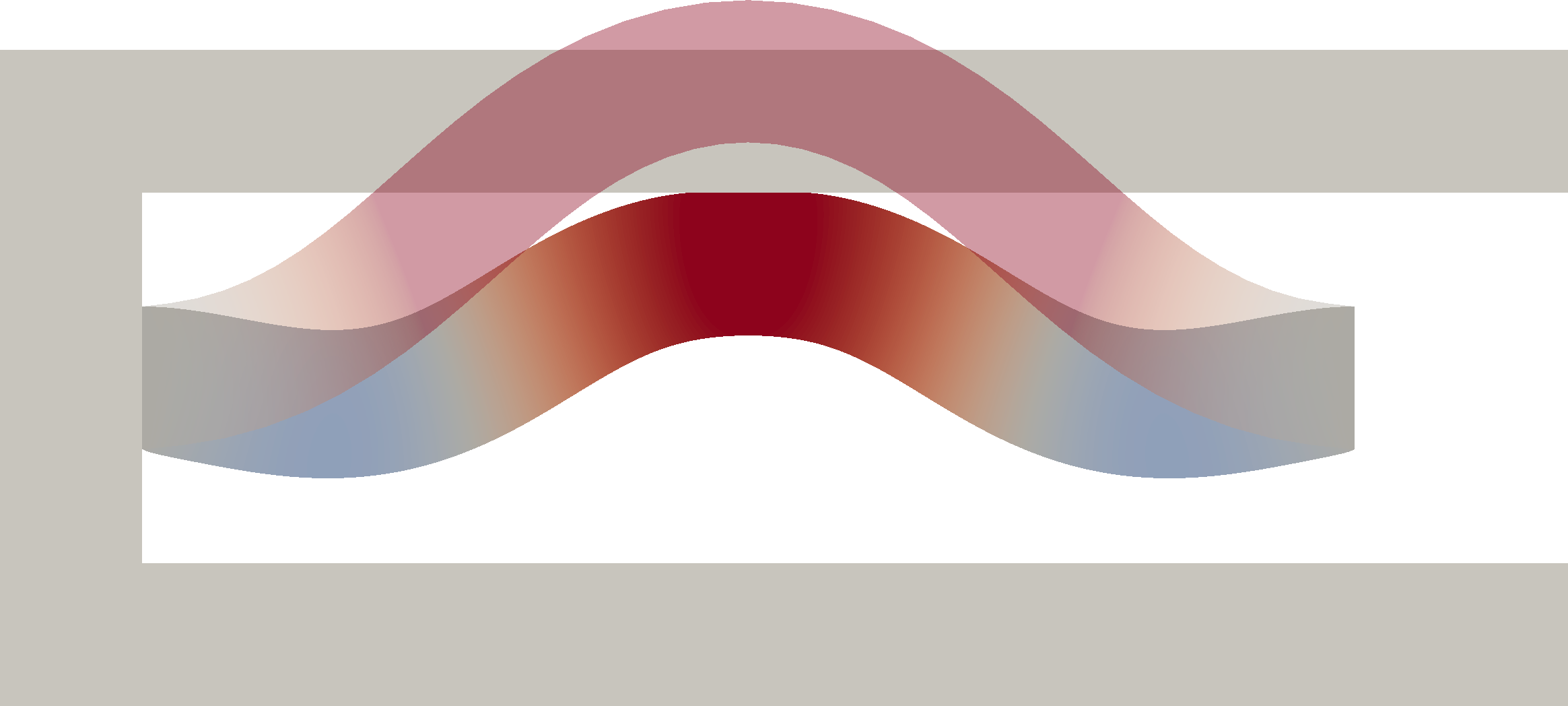} \\
\includegraphics[width=12.0cm]{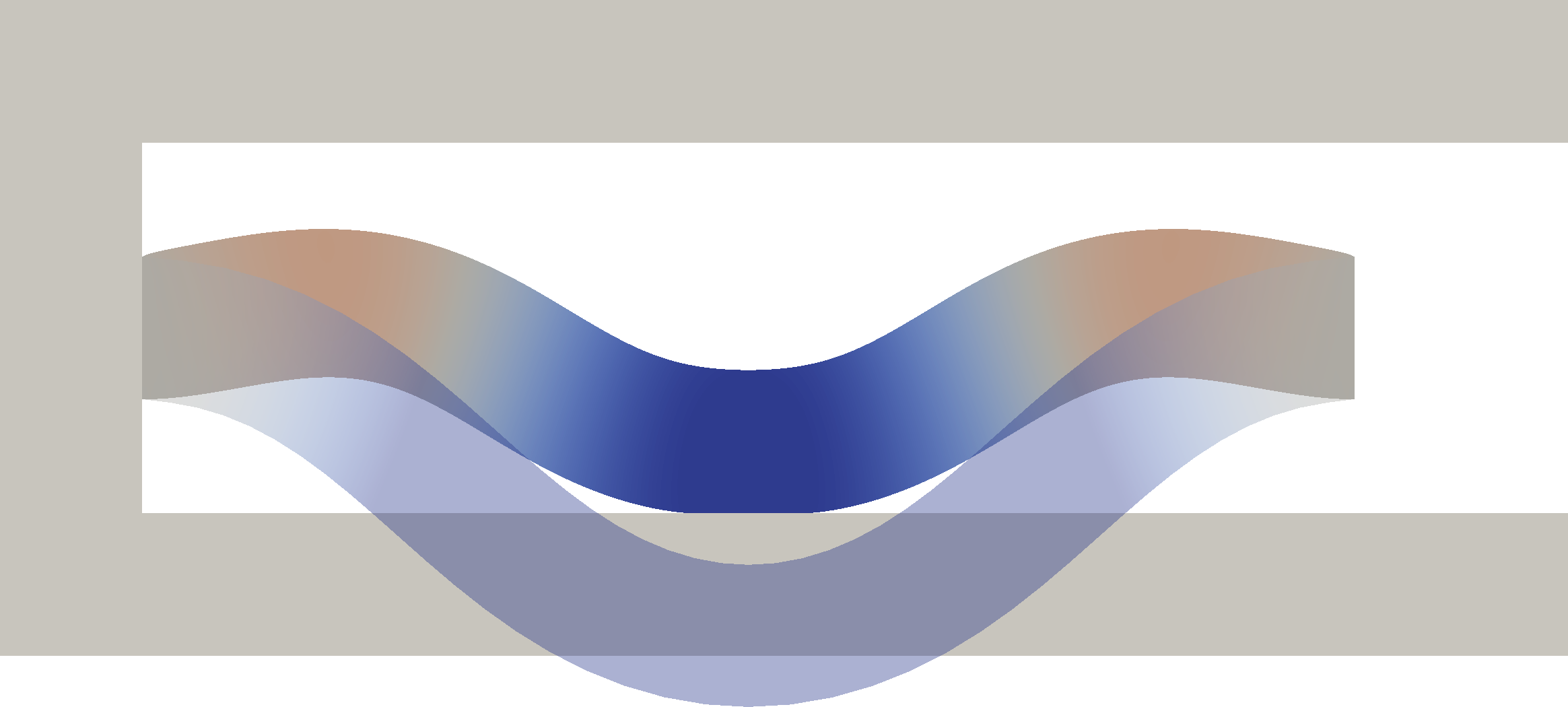} \\
\includegraphics[width=12.0cm]{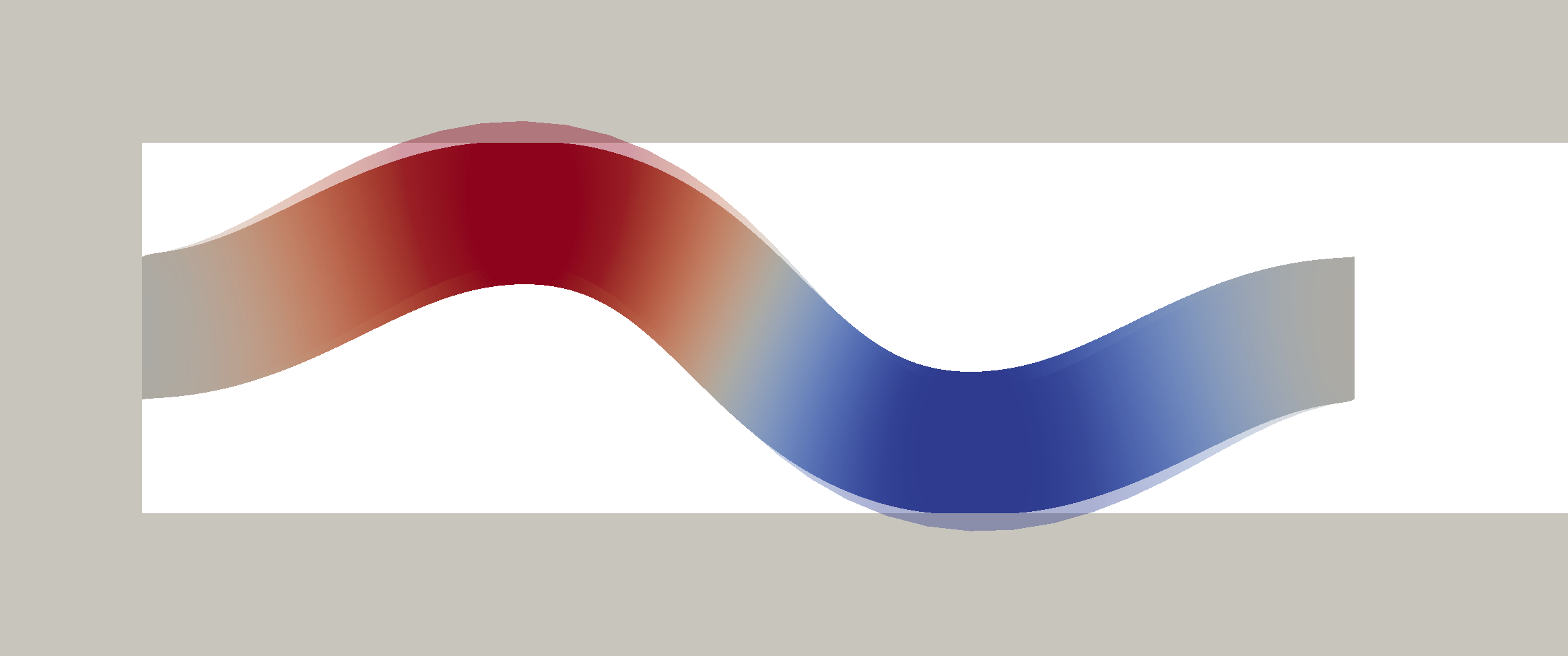} \\
\includegraphics[width=12.0cm]{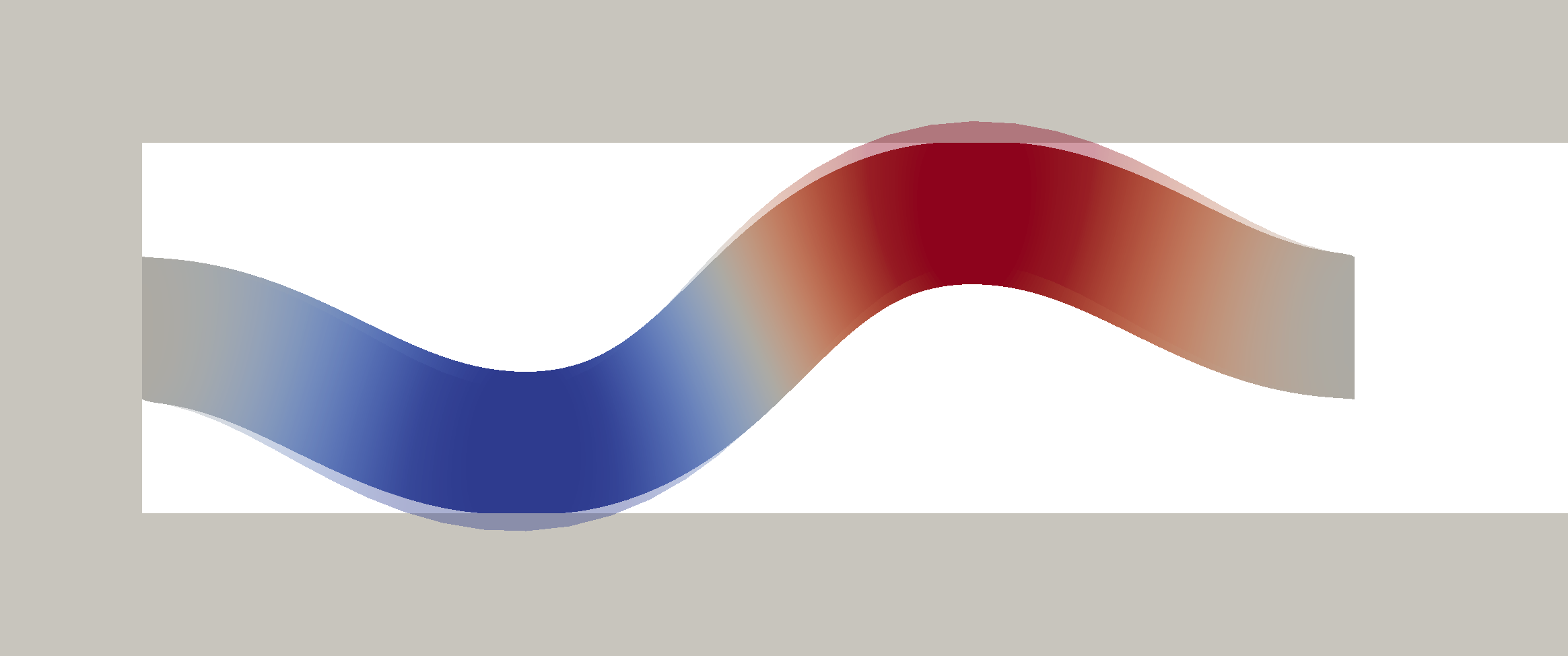} \\
\end{tabular}
\caption{Active solutions of the two-dimensional hyperelastic beam under
axial compression with obstacle constraints. The color bar denotes the vertical
component of displacement $u_2$. The corresponding \emph{unconstrained} solution
is shown semi-transparently.}
\label{fig:beam2dactive}
\end{figure}

\begin{figure}
\begin{tabular}{c}
\includegraphics[width=12.0cm]{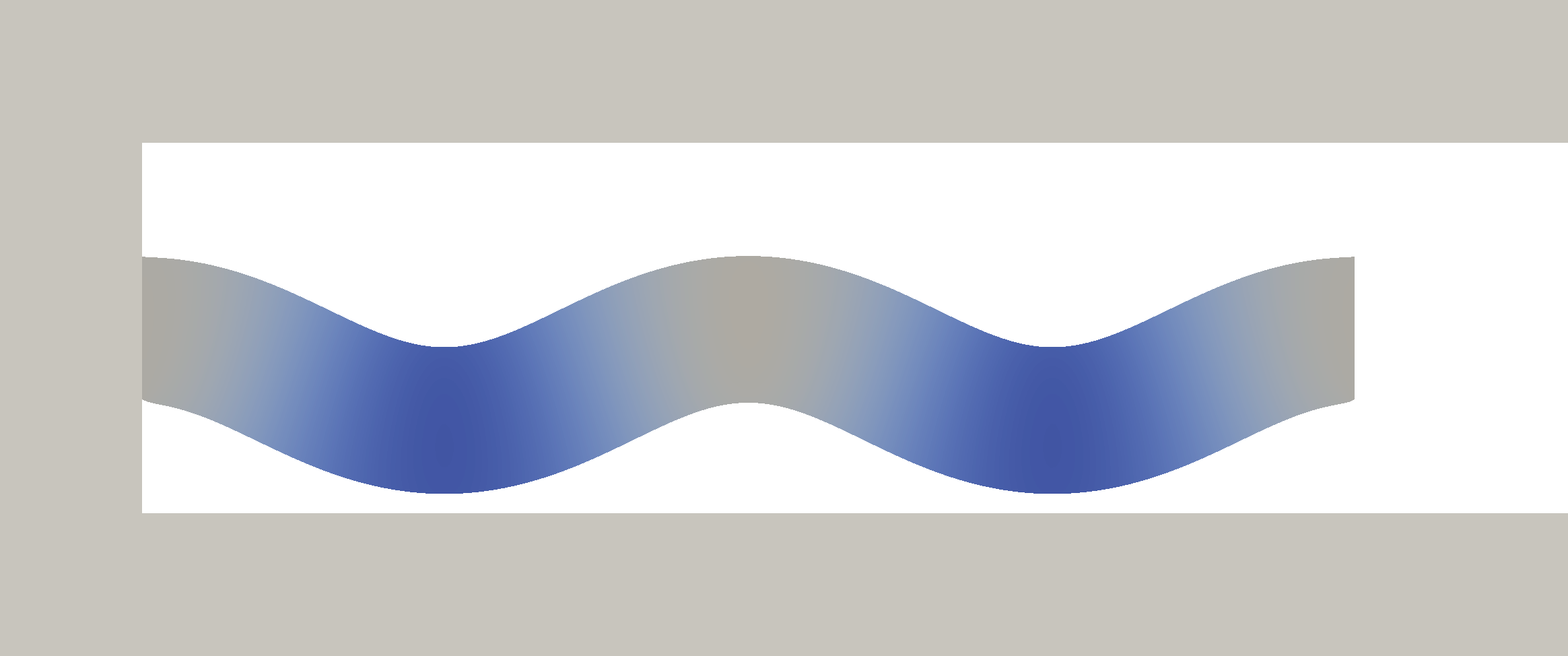} \\
\includegraphics[width=12.0cm]{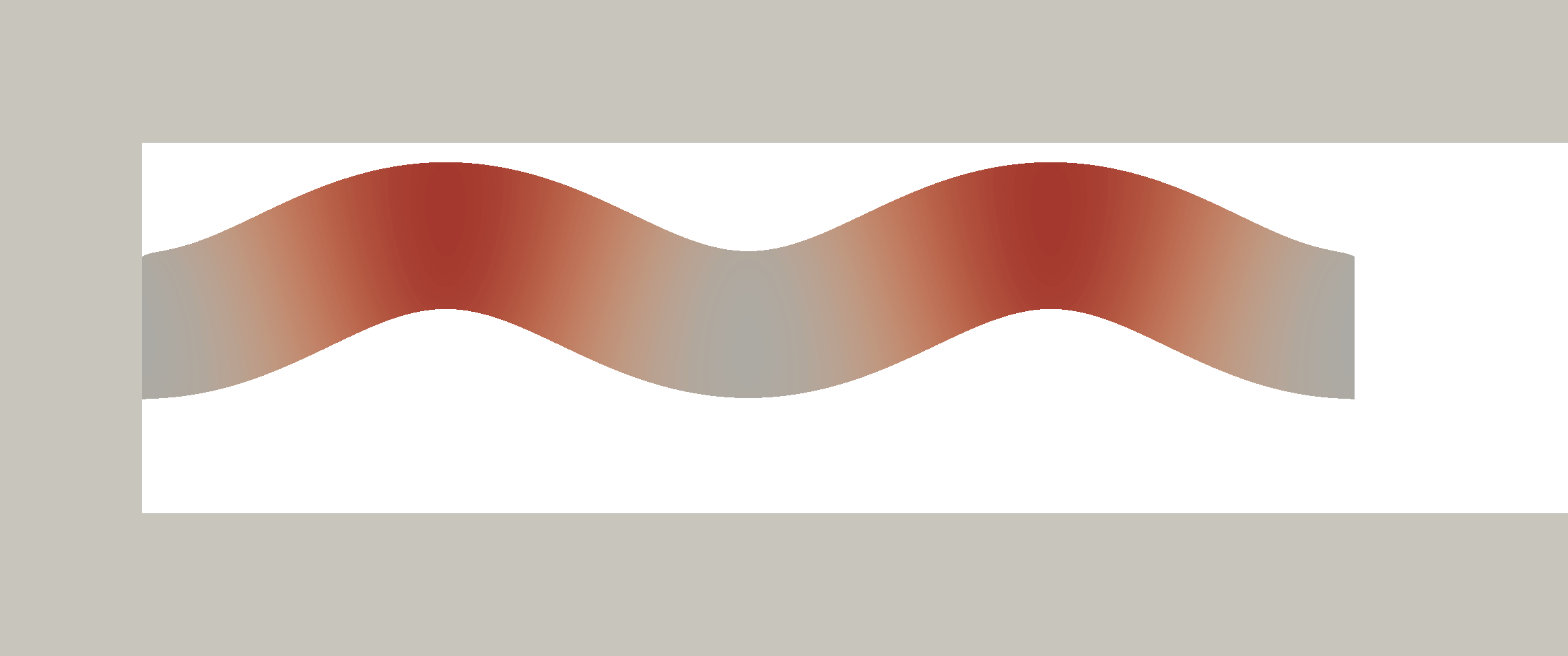} \\
\includegraphics[width=12.0cm]{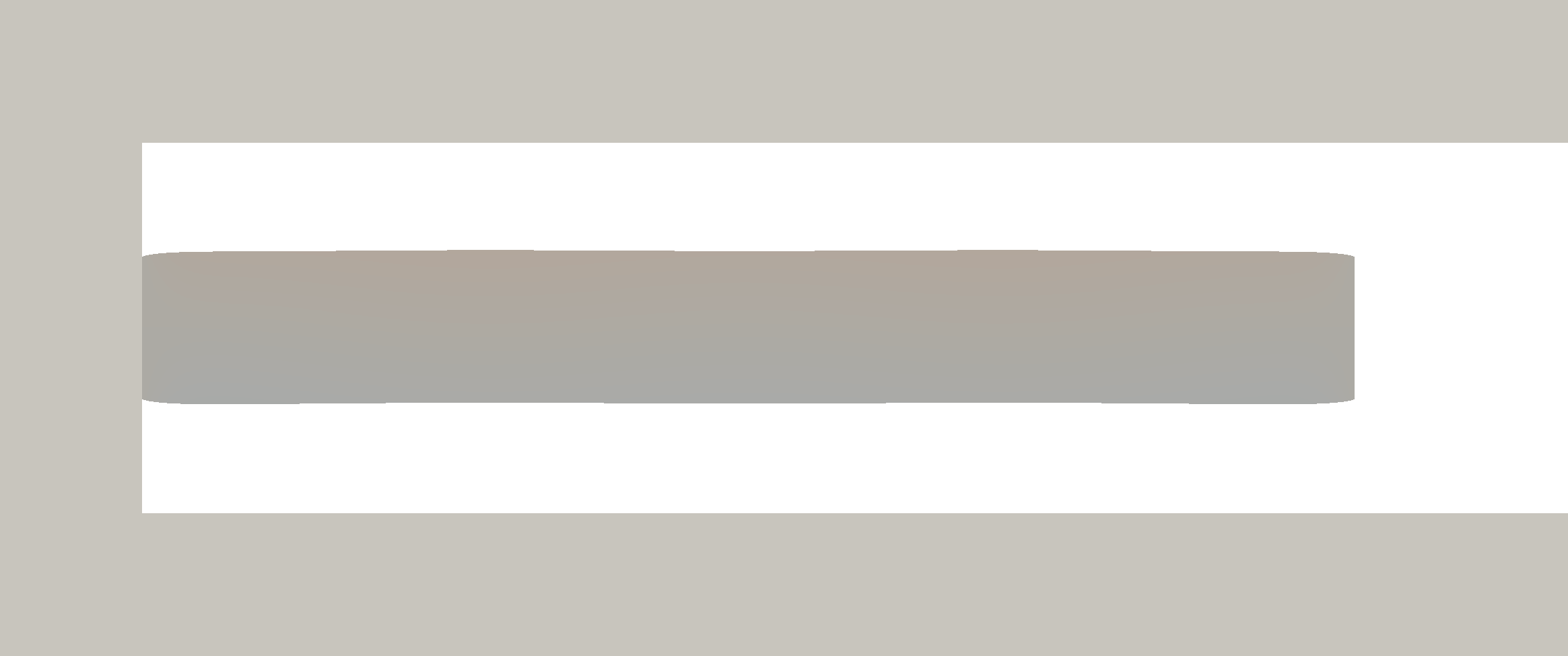} \\
\end{tabular}
\caption{Inactive solutions of the two-dimensional hyperelastic beam under
axial compression with obstacle constraints. The color bar denotes the vertical
component of displacement $u_2$.}
\label{fig:beam2dinactive}
\end{figure}

This procedure yields seven solutions for $\varepsilon = 0.15$ and $\gamma =
10^6$. To give a sense of the work involved, we report the number of semismooth
Newton iterations required for the initial discovery of each solution as
$\varepsilon$ is continued in Table \ref{tab:discoveryitercounts}, and the
average number of semismooth Newton iterations per solve and GMRES iterations per semismooth Newton step as $\gamma$ is
continued in Table \ref{tab:krylovitercounts}. The absolute and relative
tolerances of the nonlinear solver were both set to $10^{-8}$, while
the absolute and relative tolerances of the linear solver were set to
0 and $10^{-8}$. In all cases iteration counts are
modest. In particular, the results of Table \ref{tab:krylovitercounts} show
that the number of semismooth Newton iterations required does not increase as
$\gamma$ and $h$ are refined, while the number of GMRES-accelerated multigrid
V-cycles grows very slowly.

The active solutions are shown in Figure \ref{fig:beam2dactive} and the inactive
solutions are shown in Figure \ref{fig:beam2dinactive}. For each active
solution, we solve \eqref{eqn:hyper2d} \emph{without} the obstacle constraints;
the corresponding solutions are also plotted to indicate the extent to which
the obstacle constraints influence the solutions. As can be seen, the bound
constraints significantly change the solutions, and are active on a set of
positive measure on the boundary.

These results are encouraging. The function-space-based semismooth
Newton method combined with analytical path-following, parameter continuation,
multigrid and deflation appears very promising for constrained non-convex variational
problems with multiple solutions.

\section{Conclusion} \label{sec:conclusion}

Deflation is a useful technique for
identifying distinct solutions of variational inequalities
with semismooth Newton methods. In particular, employing
shifted deflation operators significantly improves the
robustness of the approach. The main strengths of the deflation method are
that it is effective, straightforward to implement and that
it does not significantly increase the cost per Newton
iteration.

While the method is found to be effective in numerical
experiments, at present no sufficient conditions are known
that guarantee convergence of the method to additional
solutions. While such conditions are unlikely to be
necessary, and may be difficult to verify \emph{a priori} in
computational practice, their availability would establish
the foundations of the method and give insight into the
design of appropriate deflation operators. The
identification of such sufficient conditions forms an
important open question and a direction for future research.

\bibliographystyle{siam}
\bibliography{\jobname.bib}

@book{nocedal2006,
  author =        {Nocedal, J. and Wright, S. J.},
  publisher =     {Springer Verlag},
  title =         {Numerical Optimization},
  year =          {2006},
}

@inbook{fichera1973,
  author =        {Fichera, G.},
  booktitle =     {Linear Theories of Elasticity and Thermoelasticity:
                   Linear and Nonlinear Theories of Rods, Plates, and
                   Shells},
  editor =        {Truesdell, C.},
  pages =         {391--424},
  publisher =     {Springer Berlin Heidelberg},
  title =         {Boundary value problems of elasticity with unilateral
                   constraints},
  year =          {1973},
  doi =           {10.1007/978-3-662-39776-3_4},
  isbn =          {978-3-662-39776-3},
}

@book{kikuchi1988,
  author =        {Kikuchi, N. and Oden, J.},
  number =        {8},
  publisher =     {SIAM},
  series =        {Studies in Applied and Numerical Mathematics},
  title =         {Contact Problems in Elasticity},
  year =          {1988},
  doi =           {10.1137/1.9781611970845},
}

@article{harker1984,
  author =        {Harker, P. T.},
  journal =       {Mathematical Programming},
  number =        {1},
  pages =         {105--111},
  title =         {A variational inequality approach for the
                   determination of oligopolistic market equilibrium},
  volume =        {30},
  year =          {1984},
  doi =           {10.1007/BF02591802},
}

@article{blowey1991,
  author =        {Blowey, J. F. and Elliott, C. M.},
  journal =       {European Journal of Applied Mathematics},
  number =        {3},
  pages =         {233--280},
  title =         {The {C}ahn--{H}illiard gradient theory for phase
                   separation with non-smooth free energy. {Part I:
                   M}athematical analysis},
  volume =        {2},
  year =          {1991},
  doi =           {10.1017/S095679250000053X},
}

@book{facchinei2003,
  author =        {F. Facchinei and J.-S. Pang},
  publisher =     {Springer-Verlag New York},
  series =        {Springer Series in Operations Research and Financial
                   Engineering},
  title =         {Finite-Dimensional Variational Inequalities and
                   Complementarity Problems},
  year =          {2003},
  doi =           {10.1007/b97543},
}

@article{ferris1997,
  author =        {Ferris, M. C. and Pang, J.-S.},
  journal =       {SIAM Review},
  number =        {4},
  pages =         {669--713},
  title =         {Engineering and economic applications of
                   complementarity problems},
  volume =        {39},
  year =          {1997},
  doi =           {10.1137/S0036144595285963},
}

@book{glowinski1984,
  author =        {R. Glowinski},
  publisher =     {Springer Berlin Heidelberg},
  title =         {Numerical methods for nonlinear variational problems},
  year =          {1984},
  doi =           {10.1007/978-3-662-12613-4},
  isbn =          {978-3-662-12615-8},
}

@book{glowinski1981,
  author =        {Glowinski, R. and Lions, J.-L. and
                   Tr\'emoli\`eres, R.},
  number =        {8},
  publisher =     {North Holland},
  series =        {Studies in Mathematics and its Applications},
  title =         {Numerical Analysis of Variational Inequalities},
  year =          {1981},
  isbn =          {9780080875293},
}

@article{harker1990,
  author =        {Harker, P. T. and Pang, J.-S.},
  journal =       {Mathematical Programming},
  number =        {1},
  pages =         {161--220},
  title =         {Finite-dimensional variational inequality and
                   nonlinear complementarity problems: A survey of
                   theory, algorithms and applications},
  volume =        {48},
  year =          {1990},
  doi =           {10.1007/BF01582255},
}

@book{kinderlehrer2000,
  author =        {Kinderlehrer, D. and Stampacchia, G.},
  number =        {31},
  publisher =     {SIAM},
  series =        {Classics in Applied Mathematics},
  title =         {An Introduction to Variational Inequalities and Their
                   Applications},
  year =          {2000},
  doi =           {10.1137/1.9780898719451},
}

@article{lions1967,
  author =        {J. L. Lions and G. Stampacchia},
  journal =       {Communications on Pure and Applied Mathematics},
  number =        {3},
  pages =         {493--519},
  title =         {Variational inequalities},
  volume =        {20},
  year =          {1967},
  doi =           {10.1002/cpa.3160200302},
}

@article{chen2000,
  author =        {X. Chen and Z. Nashed and L. Qi},
  journal =       {SIAM Journal on Numerical Analysis},
  number =        {4},
  pages =         {1200--1216},
  title =         {Smoothing methods and semismooth methods for
                   nondifferentiable operator equations},
  volume =        {38},
  year =          {2000},
  doi =           {10.1137/S0036142999356719},
}

@article{deluca1996,
  author =        {De Luca, T. and Facchinei, F. and Kanzow, C.},
  journal =       {Mathematical Programming},
  number =        {3},
  pages =         {407--439},
  title =         {A semismooth equation approach to the solution of
                   nonlinear complementarity problems},
  volume =        {75},
  year =          {1996},
  doi =           {10.1007/BF02592192},
}

@article{hintermueller2002,
  author =        {M. Hinterm\"uller and K. Ito and K. Kunisch},
  journal =       {SIAM Journal on Optimization},
  number =        {3},
  pages =         {865--888},
  title =         {The primal-dual active set strategy as a semismooth
                   {N}ewton method},
  volume =        {13},
  year =          {2002},
  doi =           {10.1137/S1052623401383558},
}

@incollection{kummer1988,
  author =        {B. Kummer},
  booktitle =     {Advances in Mathematical Optimization},
  pages =         {114--125},
  publisher =     {Akademie Verlag Berlin},
  title =         {Newton's method for non-differentiable functions},
  year =          {1988},
}

@book{klatte2002,
  author =        {D. Klatte and B. Kummer},
  publisher =     {Kluwer Academic Publishers},
  title =         {Nonsmooth Equations in Optimization},
  year =          {2002},
  doi =           {10.1007/b130810},
  isbn =          {978-1-4020-0550-3},
}

@article{qi1993,
  author =        {Qi, L. and Sun, J.},
  journal =       {Mathematical Programming},
  number =        {1},
  pages =         {353--367},
  title =         {A nonsmooth version of {N}ewton's method},
  volume =        {58},
  year =          {1993},
  doi =           {10.1007/BF01581275},
}

@book{ulbrich2011,
  author =        {Ulbrich, M.},
  publisher =     {SIAM},
  series =        {MOS-SIAM Series on Optimization},
  title =         {Semismooth {N}ewton Methods for Variational
                   Inequalities and Constrained Optimization Problems in
                   Function Spaces},
  volume =        {11},
  year =          {2011},
  doi =           {10.1137/1.9781611970692},
  isbn =          {978-1-61197-068-5},
}

@article{brown1971,
  author =        {Brown, K. M. and Gearhart, W. B.},
  journal =       {Numerische Mathematik},
  number =        {4},
  pages =         {334--342},
  publisher =     {Springer-Verlag},
  title =         {Deflation techniques for the calculation of further
                   solutions of a nonlinear system},
  volume =        {16},
  year =          {1971},
  doi =           {10.1007/BF02165004},
}

@inproceedings{levy1985,
  author =        {A. V. Levy and S. Gomez},
  booktitle =     {Numerical Optimization 1984},
  editor =        {P. T. Boggs},
  month =         {June},
  publisher =     {SIAM},
  title =         {The tunneling method applied to global optimization},
  year =          {1985},
  isbn =          {978-0898710540},
}

@phdthesis{birkisson2013,
  author =        {\'A. Birkisson},
  school =        {University of Oxford},
  title =         {Numerical Solution of Nonlinear Boundary Value
                   Problems for Ordinary Differential Equations in the
                   Continuous Framework},
  year =          {2014},
}

@article{farrell2014,
  author =        {P. E. Farrell and \'A. Birkisson and S. W. Funke},
  journal =       {SIAM Journal on Scientific Computing},
  number =        {4},
  pages =         {A2026--A2045},
  title =         {Deflation techniques for finding distinct solutions
                   of nonlinear partial differential equations},
  volume =        {37},
  year =          {2015},
  doi =           {10.1137/140984798},
}

@article{vonstengel2002,
  author =        {von Stengel, B. and van den Elzen, A. and Talman, D.},
  journal =       {Econometrica},
  number =        {2},
  pages =         {693--715},
  title =         {Computing normal form perfect equilibria for
                   extensive two-person games},
  volume =        {70},
  year =          {2002},
  doi =           {10.1111/1468-0262.00300},
}

@article{tinloi2003,
  author =        {F. Tin-Loi and P. Tseng},
  journal =       {Computer Methods in Applied Mechanics and
                   Engineering},
  number =        {11-12},
  pages =         {1377--1388},
  title =         {Efficient computation of multiple solutions in
                   quasibrittle fracture analysis},
  volume =        {192},
  year =          {2003},
  doi =           {10.1016/S0045-7825(02)00645-X},
}

@article{judice1988,
  author =        {J. J. Judice and G. Mitra},
  journal =       {European Journal of Operational Research},
  number =        {1},
  pages =         {122--128},
  title =         {An enumerative method for the solution of linear
                   complementarity problems},
  volume =        {36},
  year =          {1988},
  doi =           {10.1016/0377-2217(88)90014-8},
}

@article{conrad1988,
  author =        {F. Conrad and R. Herbin and H. D. Mittelmann},
  journal =       {SIAM Journal on Numerical Analysis},
  number =        {6},
  pages =         {1409--1431},
  title =         {Approximation of Obstacle Problems by Continuation
                   Methods},
  volume =        {25},
  year =          {1988},
  doi =           {10.1137/0725082},
}

@article{mittelmann1983,
  author =        {H. D. Mittelmann},
  journal =       {Mathematics of Computation},
  number =        {164},
  pages =         {473--485},
  title =         {An Efficient Algorithm for Bifurcation Problems of
                   Variational Inequalities},
  volume =        {41},
  year =          {1983},
}

@article{kanzow2000,
  author =        {Kanzow, C.},
  journal =       {Journal of Global Optimization},
  number =        {1},
  pages =         {1--21},
  title =         {Global optimization techniques for mixed
                   complementarity problems},
  volume =        {16},
  year =          {2000},
  doi =           {10.1023/A:1008331803982},
}

@book{allgower2003,
  author =        {Allgower, E. and Georg, K.},
  publisher =     {SIAM},
  title =         {Introduction to {Numerical Continuation Methods}},
  year =          {2003},
  doi =           {10.1137/1.9780898719154},
}

@article{fry2002,
  author =        {R. Fry and S. {McManus}},
  journal =       {Expositiones Mathematicae},
  number =        {2},
  pages =         {143--183},
  title =         {Smooth bump functions and the geometry of {B}anach
                   spaces},
  volume =        {20},
  year =          {2002},
  doi =           {10.1016/s0723-0869(02)80017-2},
}

@book{wilkinson1963,
  author =        {Wilkinson, J. H.},
  publisher =     {H.M.S.O.},
  series =        {Notes on Applied Science},
  title =         {Rounding errors in algebraic processes},
  volume =        {32},
  year =          {1963},
}

@article{mifflin1977,
  author =        {R. Mifflin},
  journal =       {{SIAM} Journal on Control and Optimization},
  number =        {6},
  pages =         {959--972},
  title =         {Semismooth and Semiconvex Functions in Constrained
                   Optimization},
  volume =        {15},
  year =          {1977},
  doi =           {10.1137/0315061},
}

@article{ralph1994,
  author =        {D. Ralph},
  journal =       {Mathematics of Operations Research},
  number =        {2},
  pages =         {352--389},
  title =         {Global Convergence of Damped
                   {N}ewton{\textquotesingle}s Method for Nonsmooth
                   Equations via the Path Search},
  volume =        {19},
  year =          {1994},
  doi =           {10.1287/moor.19.2.352},
}

@misc{farrell2015d,
  author =        {P. E. Farrell and C. H. L. Beentjes and
                   \'A. Birkisson},
  note =          {arXiv:1603.00809 [math.NA]},
  title =         {The computation of disconnected bifurcation diagrams},
  year =          {2016},
}

@article{rall1974,
  author =        {L. B. Rall},
  journal =       {SIAM Journal on Numerical Analysis},
  pages =         {34--36},
  title =         {A note on the convergence of {N}ewton's method},
  volume =        {11},
  year =          {1974},
}

@incollection{rheinboldt1978,
  author =        {W. C. Rheinboldt},
  booktitle =     {Mathematical Models and Numerical Methods},
  pages =         {129--142},
  publisher =     {Banach Center Publications},
  title =         {An adaptive continuation process for solving systems
                   of nonlinear equations},
  volume =        {3},
  year =          {1978},
}

@article{hintermueller2006,
  author =        {M. Hinterm\"uller and K. Kunisch},
  journal =       {{SIAM} Journal on Control and Optimization},
  number =        {4},
  pages =         {1198--1221},
  title =         {Feasible and Noninterior Path-Following in
                   Constrained Minimization with Low Multiplier
                   Regularity},
  volume =        {45},
  year =          {2006},
  doi =           {10.1137/050637480},
}

@article{fischer1992,
  author =        {A. Fischer},
  journal =       {Optimization},
  number =        {3--4},
  pages =         {269--284},
  title =         {A special {N}ewton-type optimization method},
  volume =        {24},
  year =          {1992},
  doi =           {10.1080/02331939208843795},
}

@inbook{kummer1992,
  author =        {Kummer, B.},
  booktitle =     {Advances in Optimization: Proceedings of the 6th
                   French-German Colloquium on Optimization Held at
                   Lambrecht, FRG, June 2--8, 1991},
  editor =        {Oettli, W. and Pallaschke, D.},
  pages =         {171--194},
  publisher =     {Springer Berlin Heidelberg},
  title =         {Newton's method based on generalized derivatives for
                   nonsmooth functions: convergence analysis},
  year =          {1992},
  doi =           {10.1007/978-3-642-51682-5_12},
  isbn =          {978-3-642-51682-5},
}

@article{kojima1986,
  author =        {M. Kojima and S. Shindo},
  journal =       {Journal of Operations Research Society of Japan},
  pages =         {352--374},
  title =         {Extensions of {N}ewton and quasi-{N}ewton methods to
                   systems of {PC}$^1$ equations},
  volume =        {29},
  year =          {1986},
}

@article{dirske1995,
  author =        {Dirkse, S. P. and Ferris, M. C.},
  journal =       {Optimization Methods and Software},
  number =        {4},
  pages =         {319--345},
  title =         {{MCPLIB}: a collection of nonlinear mixed
                   complementarity problems},
  volume =        {5},
  year =          {1995},
  doi =           {10.1080/10556789508805619},
}

@techreport{josephy1979,
  author =        {Josephy, N. H.},
  institution =   {Mathematics Research Centre, University of
                   Madison--Wisconsin},
  number =        {1965},
  title =         {{N}ewton's Method for Generalized Equations.},
  year =          {1979},
}

@unpublished{gould2002,
  author =        {N. I. M. Gould and P. L. Toint},
  note =          {19$^{\mathrm{th}}$ Biennial Conference on Numerical
                   Analysis, Dundee, Scotland},
  title =         {The state-of-the-art in numerical methods for
                   quadratic programming},
  year =          {2002},
  url =           {ftp://ftp.numerical.rl.ac.uk/pub/talks/
                  nimg.dundee.28VI01.ps.gz},
}

@article{vonneumann1945,
  author =        {von Neumann, J. and Morgenstern, O.},
  journal =       {Bulletin of the American Mathematical Society},
  number =        {7},
  pages =         {498--504},
  title =         {Theory of games and economic behavior},
  volume =        {51},
  year =          {1945},
}

@article{nash1951,
  author =        {Nash, J.},
  journal =       {Annals of Mathematics},
  number =        {2},
  pages =         {286--295},
  title =         {Non-cooperative games},
  volume =        {54},
  year =          {1951},
}

@article{lemke1964,
  author =        {Lemke, C. E. and Howson, J. T.},
  journal =       {Journal of the Society for Industrial and Applied
                   Mathematics},
  number =        {2},
  pages =         {413--423},
  title =         {Equilibrium points of bimatrix games},
  volume =        {12},
  year =          {1964},
}

@article{aggarwal1973,
  author =        {Aggarwal, V.},
  journal =       {Mathematical Programming},
  number =        {1},
  pages =         {233--234},
  title =         {On the generation of all equilibrium points for
                   bimatrix games through the {Lemke--Howson} Algorithm},
  volume =        {4},
  year =          {1973},
  doi =           {10.1007/BF01584663},
}

@book{murty1988,
  author =        {K. G. Murty},
  publisher =     {Heldermann Verlag},
  series =        {Sigma Series in Applied Mathematics},
  title =         {Linear Complementarity, Linear and Nonlinear
                   Programming},
  volume =        {3},
  year =          {1988},
  isbn =          {9783885384038},
}

@misc{gerard2017,
  author =        {H. G\'erard and V. Lecl\`ere and A. Philpott},
  note =          {arXiv:1706.08398 [math.OC]},
  title =         {On risk averse competitive equilibrium},
  year =          {2017},
}

@article{dirske1995b,
  author =        {S. P. {Dirkse} and M. C. {Ferris}},
  journal =       {Optimization Methods and Software},
  pages =         {123--156},
  title =         {The {PATH} Solver: {A} Non-Monotone Stabilization
                   Scheme for Mixed Complementarity Problems},
  volume =        {5},
  year =          {1995},
  url =           {http://www.cs.wisc.edu/~ferris/techreports/cstr1179.pdf},
}

@article{ferris2000,
  author =        {M. C. Ferris and T. S. Munson},
  journal =       {Journal of Economic Dynamics and Control},
  number =        {2},
  pages =         {165--188},
  title =         {Complementarity problems in {GAMS} and the {PATH}
                   solver},
  volume =        {24},
  year =          {2000},
  doi =           {10.1016/s0165-1889(98)00092-x},
}

@article{brune2015,
  author =        {P. R. Brune and M. G. Knepley and B. F. Smith and
                   X. Tu},
  journal =       {SIAM Review},
  number =        {4},
  pages =         {535--565},
  title =         {Composing scalable nonlinear algebraic solvers},
  volume =        {57},
  year =          {2015},
  doi =           {10.1137/130936725},
}

@book{bonnans2000,
  author =        {J. F. Bonnans and A. Shapiro},
  publisher =     {Springer New York},
  title =         {Perturbation Analysis of Optimization Problems},
  year =          {2000},
  doi =           {10.1007/978-1-4612-1394-9},
  url =           {https://doi.org/10.1007%2F978-1-4612-1394-9},
}

@article{hintermueller2004,
  author =        {M. Hinterm\"uller and M. Ulbrich},
  journal =       {Mathematical Programming},
  number =        {1},
  title =         {A mesh-independence result for semismooth {N}ewton
                   methods},
  volume =        {101},
  year =          {2004},
  doi =           {10.1007/s10107-004-0540-9},
}

@article{MHintermueller_KKunisch_2006,
  author =        {Hinterm\"uller, M. and Kunisch, K.},
  journal =       {SIAM J. Optim.},
  number =        {1},
  pages =         {159-187},
  title =         {Path-following methods for a class of constrained
                   minimization problems in function space},
  volume =        {17},
  year =          {2006},
  doi =           {10.1137/040611598},
}

@article{LAdam_HHintermueller_TMSurowiec_2018,
  author =        {Adam, L and Hinterm\"uller, M and Surowiec, T M},
  journal =       {IMA Journal of Numerical Analysis},
  title =         {A semismooth {N}ewton method with analytical
                   path-following for the {$H^1$}-projection onto the
                   {G}ibbs simplex},
  year =          {2018},
  doi =           {10.1093/imanum/dry034},
  url =           {http://dx.doi.org/10.1093/imanum/dry034},
}

@article{hintermueller2009,
  author =        {M. Hinterm\"uller and M. Hinze},
  journal =       {{SIAM} Journal on Numerical Analysis},
  number =        {3},
  pages =         {1666--1683},
  title =         {{M}oreau{\textendash}{Y}osida regularization in state
                   constrained elliptic control problems: error
                   estimates and parameter adjustment},
  volume =        {47},
  year =          {2009},
  doi =           {10.1137/080718735},
}

@book{zeidler1988,
  author =        {E. Zeidler},
  publisher =     {Springer New York},
  title =         {Nonlinear Functional Analysis and its Applications
                   IV: Applications to Mathematical Physics},
  year =          {1988},
  doi =           {10.1007/978-1-4612-4566-7},
}

@article{phu1987,
  author =        {H. X. Ph{\'{u}}},
  journal =       {Zeitschrift f{\"u}r Analysis und ihre Anwendungen},
  number =        {4},
  pages =         {371--384},
  title =         {Zur {L}\"osung des linearisierten
                   {K}nickstab-Problems mit beschr{\"a}nkter
                   {A}usbiegung},
  volume =        {6},
  year =          {1987},
  doi =           {10.4171/zaa/258},
}

@article{kirby2018a,
  author =        {R. C. Kirby},
  journal =       {{SMAI} Journal of Computational Mathematics},
  pages =         {197--224},
  title =         {A general approach to transforming finite elements},
  volume =        {4},
  year =          {2018},
  doi =           {10.5802/smai-jcm.33},
}

@misc{kirby2018b,
  author =        {R. C. Kirby and L. Mitchell},
  note =          {arXiv:1808.05513},
  title =         {Code generation for generally mapped finite elements},
  year =          {2018},
}

@article{rathgeber2016,
  author =        {F. Rathgeber and D. A. Ham and L. Mitchell and
                   M. Lange and F. Luporini and A. T. T. Mcrae and
                   G.-T. Bercea and G. R. Markall and P. H. J. Kelly},
  journal =       {{ACM} Transactions on Mathematical Software},
  number =        {3},
  pages =         {1--27},
  title =         {Firedrake: automating the finite element method by
                   composing abstractions},
  volume =        {43},
  year =          {2016},
  doi =           {10.1145/2998441},
}

@techreport{balay2017,
  author =        {S. Balay and S. Abhyankar and M. F. Adams and
                   J. Brown and P. Brune and K. Buschelman and L. Dalcin and
                   V. Eijkhout and W. D. Gropp and D. Kaushik and
                   M. G. Knepley and L. Curfman McInnes and K. Rupp and
                   B. F. Smith and S. Zampini and H. Zhang and H. Zhang},
  institution =   {Argonne National Laboratory},
  number =        {ANL-95/11 - Revision 3.8},
  title =         {{PETS}c Users Manual},
  year =          {2017},
  url =           {http://www.mcs.anl.gov/petsc},
}

@article{amestoy2001,
  author =        {P. R. Amestoy and I. S. Duff and J. Koster and
                   J.-Y. L'Excellent},
  journal =       {SIAM Journal on Matrix Analysis and Applications},
  number =        {1},
  pages =         {15--41},
  title =         {A fully asynchronous multifrontal solver using
                   distributed dynamic scheduling},
  volume =        {23},
  year =          {2001},
}

@article{ulbrich2017,
  author =        {M. Ulbrich and S. Ulbrich and D. Bratzke},
  journal =       {Journal of Computational Mathematics},
  number =        {4},
  pages =         {486--528},
  title =         {A multigrid semismooth {N}ewton method for semilinear
                   contact problems},
  volume =        {35},
  year =          {2017},
  doi =           {10.4208/jcm.1702-m2016-0679},
}

\end{document}